\documentclass[11pt]{article}
\usepackage{amsmath,amsfonts,amsthm,epsf,stackrel}
\usepackage{graphicx,subfigure}
\usepackage{natbib}
\usepackage{setspace}
\makeatletter

\usepackage{float}

\RequirePackage{natbib}

\usepackage[margin=2cm]{geometry}
\usepackage[xcolor=dvipsnames]{xcolor}

\usepackage{dsfont}

\newcommand{\ignore}[1]{}

\newcommand{\ve}{\varepsilon}

\newcommand{\x}{{\bf x}}
\newcommand{\bs}{\boldsymbol}

\newtheorem{theorem}{Theorem}[section]

\newtheorem{proposition}[theorem]{Proposition}
\newtheorem{lemma}[theorem]{Lemma}

\begin{document}
	\title{Optimal designs for the development of personalized treatment rules
		\author{David Azriel$^a$, Yosef Rinott$^b$ and Martin Posch$^c$\thanks{Corresponding author: Martin Posch; email: martin.posch@meduniwien.ac.at.}}} \maketitle
	{
		\noindent $^a$ Faculty of Industrial Engineering and Management, The Technion\\
		$^b$ Department of Statistics and Federmann Center for the Study of Rationality, The Hebrew University\\
		$^c$ Center for Medical Statistics, Informatics and Intelligent Systems, Medical University of Vienna\\
	}

	\begin{abstract}
		{We study the design of multi-armed parallel group clinical trials to estimate personalized treatment rules that identify the best treatment for a given patient with given covariates. Assuming that the outcomes in each treatment arm are given by a homoscedastic linear model, with possibly different variances between treatment arms, and that the trial subjects form a random sample from an unselected  overall population, we optimize the (possibly randomized) treatment allocation allowing the allocation rates to depend on the covariates. We find that, for the case of two treatments, the approximately optimal allocation rule does not depend on the value of the covariates but only on the variances of the responses. In contrast, for the case of three treatments or more, the optimal treatment allocation does depend on the values of the covariates as well as the true regression coefficients. The methods are illustrated with a recently published dietary clinical trial.}
	\end{abstract}
	
	\section{Introduction}

	We consider multi-armed clinical trials for selecting the best among several treatments, where the primary  outcome may depend on patients' covariates. 
	Modeling the outcome as a  function of the covariates does not only reduce the
	error variance, but allows one also to {compute and compare} ``personalized" treatment rules, that is,  rules that assign the treatment
	with the best estimated outcome for a patient having given values of the covariates.
	This work concerns the study of optimal design of clinical trials with the goal of finding an optimal personalized treatment rule.
	{We assume that a sample of subjects from the relevant population is given, and we look for an optimal allocation of subjects participating in the trial to treatments as a function of their covariates, in order to} 
	{minimize the regret for future patients. The regret is defined as the expected difference of  outcomes under the optimal and the estimated treatment rule,}
	and different designs are compared by comparing the regrets.

	{The development of statistical methods for} personalized or precision medicine is a very active area in statistics. See, e.g., \citet{KM} and \citet{KLSZ}, both containing many recent references. {Much of the literature deals with finding,  on the basis of a given sample, the best treatment as a function of personal {baseline} covariates, which {can be, for example, genetic or other biomarkers.}  Another common theme is adaptive or dynamic schemes, where treatment depends on {the individual patient's course of disease.}} 
	{\citet{KL_2019} describe two {different} settings: single- and multi-decision {problems}. Our focus here is on the former, where the treatment is only {chosen at a single time point}, and the {choice} may depend on personal {baseline} covariates. This setting is different from a classical clinical trial, in which the goal is to determine whether a certain treatment is better than control. In our setting the goal is {to propose an (approximately) optimal design of a  trial aimed at finding the best personalized treatment rule, where optimality is  with respect to the notion of the regret as defined below in \eqref{eq:regret}.}
		
		In this context, \citet{Minsker2016} and \citet{Zhao2020} study  clinical trials where subjects can be selected in an optimal way in order to obtain an efficient design. In contrast, we assume that subjects arrive according to the population distribution and cannot be selected.
		\citet{Lee2019} study optimal design for the purpose of comparing treatments, where treatment effects are estimated, and the design aims at minimizing the variances of the estimators. The covariates are assumed to be binary, and hence the design has a discrete nature.
		\citet{Zhang2021} study a similar regression model to the one considered here with two treatments and aim at finding the minimax design with respect to the variance of the decision rule. 
		Both \citet{Zhang2021} and \citet{Lee2019} differ from the present paper in the target function,  formulated in the present paper as regret, which seems to be most natural, as well as the general approach: theirs is more computational, while we aim at finding explicit optimal designs, when possible.} 
	
	In the related context of bandit problems with covariates, \citet{Golden} studied a sequential allocation scheme where at each stage one out of two treatments needs to be assigned under a minimax framework. A high-dimensional version of this problem was studied in  \citet{Bastani}. The goal in the bandit formulation is
	to determine the best treatment while  minimizing some loss or regret function for subjects in the trial sample. This is a different setting than the one considered here, since rather than be concerned with optimal treatment for subjects in the  trial, we look for designs that will allocate patients to treatments efficiently for the goal of finding the best treatment for future patients who may require one of the treatments under study.  

	Assuming that the outcomes in each treatment arm  are given by a homoscedastic linear model  (with possibly different variances between treatment arms) we  study the regret via an approximation that assumes ``ideal" conditions, such as normality instead of asymptotic normality of linear regression coefficient estimators. We show that this yields a useful approximation for quite small sample sizes. We also study an asymptotic regret. However, it turns out that when there are three or more treatments, the optimal designs for moderate sample sizes may substantially differ from the asymptotically optimal design.
	
	For the case of two treatments, we show that the optimal trial design allocates patients in proportion to the standard deviations of the response under each treatment, and  the optimal allocation probabilities do not depend on the covariates. This is shown for the case where the outcome in each treatment group depends on the covariates according to a $p$-dimensional linear regression model, as well as for the case of a single covariate and a polynomial regression model. We also formulate a minimax criterion and show that it is satisfied by designs {with allocation probabilities that do not depend on the covariates values. For parameter spaces with the same bound on the variances of the two treatments, a balanced design (equal probability) is optimal under the minimax criterion.}  
	
	For the case of three or more treatments and a linear regression model in a single covariate, we combine theoretical arguments and numerical calculations to study some examples of optimal designs. 
	We show that the optimal allocation in general depends on the covariate. For specific cases, we show that the asymptotically optimal allocation is deterministic, i.e., the range of the covariate is partitioned into intervals and in each interval the allocation probability to a particular treatment is 1. Such allocation rules arise typically in very large studies.

	{We illustrate the optimal  design in the setting of a dietary clinical trial where
		three diets  were compared in a parallel group design \citep{Ebbeling2018}.
		The primary outcome in the original trial was averaged total energy expenditure, which is correlated with weight loss. An important patient {covariate} is {insulin secretion} and we consider  the objective  of determining  a treatment rule that assigns for each patient the best diet based on his or her insulin secretion level measured at baseline. Here the utility of the diets is defined by the  averaged total energy expenditure adjusted by diet costs.
		We derive the optimal experimental design, that is, the optimal allocation of subjects to treatments as a function of the measured insulin level, minimizing a suitable regret function.  
		Based on the data in \cite{Ebbeling2018} we quantify the improvement in regret of the optimal design compared to designs with optimal allocation rates that do not depend on the covariate and the  design with equal allocation rates that was used in the original trial (see Section \ref{sec:diets} for details).}

	{The paper is organized as follows. In Section \ref{sec:problem} the design problem is formalized and in particular the regret defined.   In Section \ref{sec:compreg} an explicit formula that approximates the regret (dubbed the \textit{ideal regret}) is derived. The optimal design for the case of two treatments is discussed in Section \ref{sec:two_treatments} after studying the rate of the approximating ideal regret to the true one. We also compute the  asymptotic regret, give a numerical example, and discuss the case of polynomial regression (Section \ref{sec:polynomial}). In Section \ref{sec:K_treatments} the setting of $K$ treatments and a single covariate is considered and the corresponding asymptotic   regret is derived. In Section \ref{sec:Knorm} we derive formula \eqref{eq:regretdirectdetailmdy} for the ideal regret, and we compare the regret and ideal regret asymptotically. Section \ref{sec:numopt} discusses an algorithm based on \eqref{eq:regretdirectdetailmdy} to obtain an approximate optimal design, and demonstrates the reduction in regret by such optimization.  Furthermore, in Section \ref{sec:lowerbound} the computation of a lower bound for the optimal  regret is discussed and an example is given where  the optimal allocation rule minimizing the asymptotic regret can be derived. In Section \ref{sec:diets} the procedure is illustrated with an example, and in Section \ref{sec:discussion} limitations and possible extensions of the approach are discussed.}

	
	\section{Problem Statement \label{sec:problem}}
	Consider  $K$ possible treatments $T_1, \ldots, T_K$. Let $Y$ be a one-dimensional continuous response variable and let ${\bf X} \in \mathbb{R}^p$ be a vector of a subject's covariates. We assume that the joint distribution of the covariates is continuous with density denoted by $f(\x)$. {If some of the covariates are discrete,  the density $f(\x)$ should be taken with respect to a suitable measure, and our results apply }{as long as one of the covariates has a continuous distribution (conditioned on all other covariates).} {For simplicity of notation, we assume ${\bf X}$ has a density with respect to Lebesgue measure. } 
	The expected response is $E(Y|{\bf X},T_k)=g_k({\bf X})$, where $g_k$ is an unknown function. The optimal treatment for a subject with a covariate vector ${\bf x}$ is $\delta^*({\bf x})=\arg \max_{k \in \{1,\ldots,K\}} g_k({\bf x})$ (assuming that a higher response is better). If the above arg max contains more than one $k$, one arbitrary treatment is selected. 
	
	Suppose that a clinical trial with $n$ subjects is performed in order to estimate $\delta^*$. Let ${\bf X}_1,\ldots,{\bf X}_n$ denote a sample according to $f$ of the covariates of the $n$ subjects. The design we study here consists of allocating subjects to treatments, taking account of their covariates. Each subject $i$ is allocated independently  to treatment $k$ with probability $\pi_k({\bf X}_i)$, where $\pi_1({\bf x}),\ldots,\pi_K({\bf x})$ are non-negative functions  satisfying for each ${\bf x}$, $\sum_{k=1}^K \pi_k({\bf x})=1$. The allocation functions define densities for the covariates in the $K$ treatment groups. Specifically, in treatment $k$, ${\bf X}$ is sampled from the density $f_k({\bf x}):= f({\bf x}) \pi_k({\bf x})/ \nu_k$ where $\nu_k:=\int f({\bf x}) \pi_k({\bf x}) d{\bf x}$. {Note that using the sampled ${\bf X}_1,\ldots,{\bf X}_n$ we can estimate the density $f$ using a parametric or nonparametric estimator. In fact, we may have {more  observations of the covariates} ${\bf X}$  than those $n$ values which will enter the trial.}
	Our purpose is to find the optimal allocation functions $\pi_1({\bf x}),\ldots,\pi_K({\bf x})$ that minimize the regret, which is defined next.

	Let $\widehat{\delta}$ denote the final estimate of the optimal rule $\delta^*$ based on the covariates, responses, and allocations of the subjects in the trial. Then the {\it regret} for a randomly chosen future patient is defined as
	\begin{equation} \label{eq:regret}
		R(\pi_1,\ldots,\pi_K):=E\left[ g_{\delta^*(\widetilde{\bf X})}(\widetilde{\bf X}) - g_{\widehat{\delta}(\widetilde{\bf X})}(\widetilde{\bf X})\right],
	\end{equation}
	where $\widetilde{\bf X}$ denotes the covariate vector of the chosen patient and the expectation is over $\widetilde{\bf X}$, $\widehat{\delta}$, and the randomizations in the allocations. In  words, the regret is the difference between the expected response to the optimal treatment $\delta^*$ and its estimate $\widehat{\delta}$ for any independent future subject arising from the same population as those in the trial.  
	A similar criterion appears, for examples, in  \citet{Qian2011}; see also \citet{KM} and references therein.  
	Our goal is to  minimize the regret over all feasible allocations, that is, nonnegative functions $\pi_k(\cdot)$ that satisfy $\sum_{k=1}^K \pi_k(\x)=1$ for every $\x$.
	
	In this study we assume that we have some preliminary knowledge about the parameters of the problem, such as the functions $g_k$, and the variances of $Y$ given ${\bf X},T_k$ and the density $f$ of the covariates.
	This approach concurs with existing literature on locally optimal designs; see,  e.g., \citet{Chernoff} and \citet[Chapter~6]{Silvey}, where  optimality is achieved for a given set of parameters. \citet{Sverdlov2013} review several methodological advances in (local) optimal allocation for clinical trials. Our goal is to find optimal allocation functions relative to the assumed parameter values, and obtain a treatment rule $\widehat{\delta}$ based on an experiment according to our design. Knowledge about the parameters can arise from previous experience or theory,  from a pilot studies, or in some situations from earlier phases of the study. However, we do not consider a sequential approach in which subjects arrive sequentially and each one is assigned to a treatment that will make the overall trial optimal.

	\section{The regret under ``ideal" conditions}\label{sec:compreg}
	Let $(Y_1,{\bf X}_1),\ldots, (Y_n,{\bf X}_n)$ be a sample obtained by certain allocation functions $\pi_1,\ldots,\pi_n$, given along with the treatment allocated for each ${\bf X}_i$. 
	Assume $K\geq 2$ and ${\bf X} \in  \mathbb{R}^p$, and under treatment $T_k$ we have $Y=g_k({\bf X})+\varepsilon$, where $g_k({\bf X})=E(Y|{\bf X},T_k)=\alpha_k+{\bs \beta}_k^t{\bf X}$, and 
	$\sigma_k^2:= Var(\varepsilon|{\bf X},T_k)$. 
	Let $\widehat{\alpha}_k, {\widehat{\bs \beta}_k}$ denote the OLS estimators based on data from treatment $T_k$, $\widehat g_k({\bf X})$ the corresponding estimated regression functions, and $\widehat{\delta}({\bf x}) = \arg\max_k \widehat g_k({\bf X})$ the decision rule. Then,
	\begin{align} \label{eq:regretdirect}
		R(\pi_1,\ldots,\pi_K)&=\sum_{k=1}^K\int_{\mathbb{R}^p}P\left({\widehat{\delta}({\bf x})}=k
		\right)\left[ g_{\delta^*({\bf x})}({\bf x}) - g_k({\bf x})\right]f({\bf x})d{\bf x}.
	\end{align}	
	For $k=1,\ldots, K$, let 
	\begin{equation}\label{eq:def_Q}
		{\bs Q}_k := \int_{\mathbb{R}^p} \binom{1}{\bf x} ({1},{\bf x}^t) f_k({\bf x}) d{\bf x},\quad {\bs \Sigma}_{k}
		:=\frac{1}{\nu_k}\sigma^2_k {\bs Q}_k^{-1},\quad \xi_{k}^2({\bf x}):=({1},{\bf x}^t) {\bs \Sigma}_{k} \binom{1}{\bf x}.    
	\end{equation}
	It is well known that under standard conditions (see, e.g., \citet{Hansen}  Chapter 7) we have $\sqrt{n}[(\widehat{\alpha}_k, {\widehat{\bs \beta}_k})-({\alpha}_k,{\bs \beta}_k)] \to N(0, {\bs \Sigma}_{k})$, and in particular  ${\bs \Sigma}_{k}$ is the asymptotic variance of the OLS estimators.
	
	In order to approximate the probability $P\left({\widehat{\delta}({\bf x})}=k\right)$, we first impose the ``{\em ideal}" conditions, namely that $\widehat{g}_{k}({\bf x})$ possess their joint asymptotic distribution. We shall later drop this assumption .
	In particular, by considering the regression model $Y=\sum_{k=1}^K  \{\alpha_k+{\bs \beta}_k^t{\bf X}\}{\mathds{1}_k} {+\ve}$, where  ${\mathds{1}_k}$ is the indicator of treatment $k$, we obtain that under  the ideal conditions, the OLS estimators are jointly normal with  $Var(\widehat\alpha_k,\widehat{\bs \beta}_k)={\bs \Sigma}_{k}$ and hence $Var(\widehat g_k({\bf x}))=\xi_{k}^2({\bf x})/n$. Also, under the ideal conditions the estimators $\widehat{g}_{k}({\bf x})$ are independent conditionally on  ${\bf X}_1,\ldots {\bf X}_n, T_1, \ldots T_K$, with constant expectation $\alpha_k+{\bs \beta}_k^{t}{\bf x}$. By the law of total covariance, conditional independence and constant conditional expectations imply that   $\widehat{g}_{k}({\bf x})$ are uncorrelated, and by joint normality, they are independent.

	Therefore, with $P_I$ denoting the  probability under the ideal conditions, we have
		{	   
			\begin{align}  \label{eq:probsdy}
				P_I\left({\widehat{\delta}({\bf x})}=k\right)&=P_I\left(\max_{\ell=1,\ldots,K,\ell\not =k} \widehat{g}_{\ell}({\bf x})<\widehat g_k({\bf x})\right) \nonumber  \\ \nonumber&=  P_I\left( \bigcap_{\ell \neq k} \left\{ \frac{\sqrt{n}[\widehat{g}_{\ell}({\bf x})-{g}_{\ell}({\bf x})]}{\xi_{\ell}({\bf x})}<\frac{\sqrt{n}[\widehat{g}_{k}({\bf x})-{g}_{\ell}({\bf x})]}{\xi_{\ell}({\bf x})} \right\} \right )\\
				\nonumber&=P\left( \bigcap_{\ell \neq k} \left\{ {Z}_{\ell} <\frac{Z_k\xi_{k}(\x)+ \sqrt{n}[{g}_{k}({\bf x})-{g}_{\ell}({\bf x})]}{\xi_{\ell}(\x)} \right\}\right)\\
				\nonumber&=\int P\left( \bigcap_{\ell \neq k} \left\{ {Z}_{\ell} <\frac{z\xi_{k}(\x)+ \sqrt{n}[{g}_{k}({\bf x})-{g}_{\ell}({\bf x})]}{\xi_{\ell}(\x)} \right\} \Big| Z_k=z \right)\varphi(z)dz \\
				&=\int \prod_{\ell=1,\ldots,K,\ell\not=k} \Phi\left(\frac{{z}\xi_{k}({\bf x})+\sqrt{n}[g_k({\bf x})-g_\ell({\bf x})]}{\xi_{\ell}({\bf x})}\right)\varphi(z)dz,
			\end{align}
			where $P$ is with respect to $\{Z_m\}_{m=1}^K$, which are independent $N(0,1)$, ${Z}_m$ represents $\frac{\sqrt{n}[\widehat{g}_{m}({\bf x})-{g}_{m}({\bf x})]}{\xi_m({\bf x})}$; $\varphi$ and  $\Phi$ denote the standard normal density and cumulative distribution function.
			The product in the last line was justified above by independence under the ideal conditions, which is also used to compute the conditional probability.} We define the \textit{ideal regret}, denoted by $R_I$, by
		\begin{equation}\label{eq:IR}
			R_I(\pi_1,\ldots,\pi_K):=\sum_{k=1}^K\int_{\mathbb{R}^p}P_I\left({\widehat{\delta}({\bf x})}=k
			\right)\left[ g_{\delta^*({\bf x})}({\bf x}) - g_k({\bf x})\right]f({\bf x})d{\bf x},
		\end{equation}
		where $P_I\left({\widehat{\delta}({\bf x})}=k\right)$ is given in \eqref{eq:probsdy}.

		\section{The case of two treatments}\label{sec:two_treatments}


		\subsection{Regret approximations and optimal design }\label{sec:approx_regret}
		We assume now that $K=2$ and  $Y=g_k({\bf X})+\varepsilon$ under treatment $T_k$ for $k=1,2$, where $g_k({\bf X})$ denotes the conditional mean of $Y$ given ${\bf X},T_k$, where ${\bf X} \in \mathbb{R}^p$ and 
		\begin{equation}\label{eq:model_two}
			g_1({\bf X})=\alpha_1+{\bs \beta}_1^t {\bf X},~g_2({\bf X})=\alpha_2+{\bs \beta}_2^t {\bf X},~Var(\ve|{\bf X},T_1)=\sigma_1^2~\text{and}~Var(\ve|{\bf X},T_2)=\sigma_2^2.  
		\end{equation}
		We do not assume normality of the errors. However, we assume the existence of the moment generating function of $\ve$ when conditioned on ${\bf X},T_k$. This assumption is needed for some large deviation exponential bounds used below, and could be relaxed to assuming finiteness of some high order moments instead. We also assume that ${\bf X}$ is continuous with a bounded density $f({\bf x})$  supported on $[0,1]^p$.

		Under the ideal conditions of Section \ref{sec:compreg} with  $K=2$, Equation \eqref{eq:probsdy} becomes 
		\begin{equation}\label{eq:def_V}
			P_I\left({\widehat{\delta}({\bf x})}=k\right) = \Phi\left(\frac{\sqrt{n}[g_k({\bf x})-g_l({\bf x})]}{\sqrt{V({\bf x})}}\right)\!\!,\text{where } V({\bf x}):= {\xi^2_{1}({\bf x})+\xi^2_{2}({\bf x})}=(1,{\bf x}^t)\left( {\bs \Sigma}_{1}+ {\bs \Sigma}_{2} \right)\binom{1}{{\bf x}}
		\end{equation}
		for $1\le k \neq l \le 2$.
		The  ideal regret (\ref{eq:IR})  is
		\begin{align}\label{eq:reg_approx} R_I(\pi_1,\pi_2)=&\int_{\{{\bf x}: g_1({\bf x})>g_2({\bf x})\}}\Phi\left(\frac{\sqrt{n}[g_2({\bf x})-g_1({\bf x})]}{\sqrt{V({\bf x})}}\right) [ g_{1}({\bf x}) - g_2({\bf x})]f({\bf x})d{\bf x}\nonumber\\&+\int_{\{{\bf x}: g_1({\bf x})<g_2({\bf x})\}}\Phi\left(\frac{\sqrt{n}[g_1({\bf x})-g_2({\bf x})]}{\sqrt{V({\bf x})}}\right) [ g_{2}({\bf x}) - g_1({\bf x})] f({\bf x})d{\bf x}\nonumber\\&
			=\int_{[0,1]^p}\Phi\left(-\frac{\sqrt{n}|g_1({\bf x})-g_2({\bf x})|}{\sqrt{V({\bf x})}}\right) |g_{2}({\bf x}) - g_1({\bf x})| f({\bf x})d{\bf x}.
		\end{align}
		Note that $R_I$ depends also on the parameters $\alpha_1,\alpha_2,{\bs \beta}_1,{\bs \beta}_2,\sigma^2_1,\sigma^2_2$, which are suppressed.
		By \eqref{eq:reg_approx}, minimization of $R_I(\pi_1,\pi_2)$ amounts to finding the allocation functions $\pi_1,\pi_2$ that minimize $V(\x)$. The following theorem states that this can be done by minimizing $V(\x)$ uniformly over $\x$.
		\begin{theorem}\label{thm:v}
			The allocation functions $\pi^0_1({\bf x}):=\frac{\sigma_1}{\sigma_1+\sigma_2}$ and $\pi^0_2({\bf x}):=\frac{\sigma_2}{\sigma_1+\sigma_2}$ minimize both $V(\x)$ uniformly over $\x$, and $R_I(\pi_1,\pi_2)$.    
		\end{theorem}
		{Note that in order to implement the allocation functions $\pi^0_1({\bf x})$ and $\pi^0_2({\bf x})$ one only needs to know the ratio of the variances of the errors of the two treatments, and no other parameters of the model.} {The allocation functions also does not require knowledge of the density $f(\x)$.}

		Recalling that $V(\x)$ can be written as a quadratic form,  see \eqref{eq:def_V}, Theorem \ref{thm:v} follows immediately from Lemma   \ref{lem:order}	below whose proof, as all other proofs, is given in the Appendix.	
		Given matrices ${A}$ and ${B}$,
		we write $ {A} \succeq {B}$ if ${A}-{B}$ is positive semi-definite. Let ${\bs Q} := \int_{[0,1]^p}  \binom{1}{\bf x} ({1},{\bf x}^t) f({\bf x}) d{\bf x}$, and recall the notation of \eqref{eq:def_Q}. We have
		\begin{lemma}\label{lem:order}
			With the above definitions,  
			\begin{equation}\label{eq:ineq_hetro}
				{\bs \Sigma}_{1}+ {\bs \Sigma}_{2} \succeq (\sigma_1+\sigma_2)^2 \mathbb{\bs Q}^{-1}.
			\end{equation}
			Moreover, the lower bound is attained when $\pi_1({\bf x})=\pi^0_1({\bf x})$ and $\pi_2({\bf x})=\pi^0_2({\bf x})$.
		\end{lemma}
		To see the latter statement note that when $\pi_k({\bf x})$ do not  depend on $\bf x$, both matrices ${\bs Q}_k$ are proportional to ${\bs Q}$ and therefore both matrices ${\bs \Sigma}_{k}$ are proportional to ${\bs Q}^{-1}$.
		
		Instead of assuming that the regression parameters are known, consider now a minimax criterion for the ideal regret. Let ${\bs \gamma}$ denote the regression parameters, i.e., ${\bs \gamma}:=(\alpha_1,\alpha_2,{\bs \beta}_1,{\bs \beta}_2,\sigma_1,\sigma_2)$ and consider the parameter space ${ \Gamma}$, where for some constant $S$,
		$\sigma_1,\sigma_2 \le {S}$ and the rest of the parameters are unrestricted. The bound on the variances is required to make the problem nondegenerate. The allocation functions $\pi_1^*,\pi_2^*$ are said to be minimax if 
		\begin{equation}\label{eq:minimax}
			\sup_{\bs\gamma \in \Gamma} R_I(\pi^*_1,\pi^*_2) \le \sup_{\bs\gamma \in \Gamma} R_I(\pi_1,\pi_2) \text{ for all allocation functions }\pi_1,\pi_2.
		\end{equation}
		By \eqref{eq:reg_approx} and the definition of $V(\x)$ in \eqref{eq:def_V} it is easy to see that the regret is maximized for ${\bs\gamma \in \Gamma}$ when  $\sigma_1=\sigma_2={S}$, in which case Theorem \ref{thm:v} implies that the optimal design is $\pi_1^*(\x)=\pi_2^*(\x)=1/2$. This statement is true regardless of the value of the other parameters in $\bs \gamma$. It follows that the balanced design satisfies \eqref{eq:minimax} and is therefore minimax. The same argument shows that if the bound $\sigma_1,\sigma_2 \le {S}$ is replaced by $\sigma_1\le S_1,\sigma_2\le {S}_2$, then the minimax design is $\pi_1^*(\x)=\frac{{S_1}}{{S_1} + {S_2}}, \pi_2^*(\x)=\frac{{S_2}}{{S_1} + {S_2}}$.   
		
		{If $\sigma_1$ and $\sigma_1$ are considered unknown and reliable estimates are not available then one should use the minimax balanced design. We now consider the sensitivity of the regret to misspecification of the standard deviations.  By Theorem \ref{thm:one_dim_lim} below the asymptotic regret is given by $\displaystyle	
			\lim_{n \to \infty} n R(\pi_1,\pi_2)= \frac{V(\theta) f(\theta)}{2 (\beta_2-\beta_1)}{I\{\theta \in [0,1]\}}$. The ratio of the asymptotic regret using  the balanced design and the optimal design  given standard deviations $\sigma_1$ and $\sigma_2$ and $\pi_1=\sigma_1/(\sigma_1+\sigma_2)$ is readily seen to be $$\frac{R(1/2,1/2)}{R(\pi_1,\pi_2)}=
			\Big(\frac{\sigma_1^2}{1/2}+\frac{\sigma_2^2}{1/2}\Big) \Big/ \Big(\frac{\sigma_1^2}{\pi_1}+\frac{\sigma_2^2}{\pi_2}\Big).$$
			If, for example, $\sigma_1/\sigma_2=1.5$ then the optimal design has $\pi_1=1.5/2.5=0.6$, and by the above formula 
			$\displaystyle \Big(\frac{\sigma_1^2}{1/2}+\frac{\sigma_2^2}{1/2}\Big) \Big/ \Big(\frac{\sigma_1^2}{\pi_1}+\frac{\sigma_2^2}{\pi_2}\Big)=1.04$, a loss of 4\% incurred by using the balanced rather than the optimal design.
			This shows that the asymptotic regret is not very sensitive to initial estimates or assumptions on the standard deviations, and the balanced design performs well even if the true standard deviations are quite far from being equal. If $\sigma_1/\sigma_2=2$, and we use the balanced design the ratio of the regrets is 1.11, which means that using the non-optimal balanced design rather than the optimal one with $\pi_1=2/3$ increases the regret by about 11\%.}
	{Similar calculations for the ideal regret \eqref{eq:reg_approx} are more complicated because the result depends on the parameters of the problem.  In Section \ref{sec:2_numeric_example} below a specific example is considered, and it is shown that when using the balanced design, rather than the optimal one, the ideal regret increases by 3\%. These findings indicate that using the balanced design could yield a regret that is close to optimal in reasonable scenarios.}

	We now discuss the regret \eqref{eq:regretdirect}. In Theorem \ref{thm:mult_K_2} we  approximate the regret by the ideal regret, and then we use the above results for the ideal regret to find an asymptotically optimal design for the regret itself. We later demonstrate numerically that the latter design provides a good approximation for finite sample sizes. 
	The regret is given by
	\begin{multline}\label{eq:reg2}
		R(\pi_1,\pi_2)=\int_{\{ {\bf x} : g_1({\bf x}) > g_2({\bf x})\}} P(\widehat{g}_2({\bf x}) > \widehat{g}_1({\bf x}))(g_1({\bf x}) - g_2({\bf x}))f({\bf x}) d{\bf x}\\
		+\int_{\{{\bf x} : g_2({\bf x}) > g_1({\bf x})\}} P(\widehat{g}_1({\bf x}) > \widehat{g}_2({\bf x}))(g_2({\bf x}) - g_1({\bf x}))f({\bf x}) d{\bf x}.
	\end{multline}

	\begin{theorem}\label{thm:mult_K_2}
		Under model \eqref{eq:model_two} with $\nu_1,\nu_2>0$, we have for any $\varepsilon>0$
		\begin{equation*}
			n^{3/2-\varepsilon} \left| R(\pi_1,\pi_2) - R_I(\pi_1,\pi_2)	\right| \to 0 \text{  as  } n \to \infty.
		\end{equation*}
	\end{theorem}

	To facilitate the representation, we state the next result in the one-dimensional case first.  
	In this case $g_1$ and $g_2$ are functions of a single {continuous} variable, and we assume   that their intersection point $\theta :=\frac{\alpha_1-\alpha_2}{\beta_2-\beta_1}$ is in $[0,1]$. 
	Assume without loss of generality, that $\beta_2>\beta_1$. It follows that for $x<\theta$, $g_1(x)>g_2(x)$, that is, $T_1$ is the better treatment, and when $x>\theta$, $T_2$ is better. In the case $p=1$ the limit of the regret is simple, given by the following theorem. 
	\begin{theorem}\label{thm:one_dim_lim}
		Under model \eqref{eq:model_two} with $p=1$,   if $\nu_1,\nu_2>0$, then\,\,\,\,
		$\displaystyle	
		\lim_{n \to \infty} n R(\pi_1,\pi_2)= \frac{V(\theta) f(\theta)}{2 (\beta_2-\beta_1)}{I\{\theta \in [0,1]\}}$.
	\end{theorem}
	{Notice that if $\beta_1=\beta_2$, then $\lim_{n \to \infty} n R(\pi_1,\pi_2)=0$, because in this case $R(\pi_1,\pi_2)$ decreases to zero at an exponential rate. In our setting, $\beta_1 \ne \beta_2$ and when $\theta \in [0,1]$,  $\lim_{n \to \infty} n R(\pi_1,\pi_2)$ is positive. Theorem \ref{thm:one_dim_lim} implies that in order to minimize this limit} we have to minimize $V(\theta)$, which by Theorem \ref{thm:v} is achieved  by taking $\pi_1({\bf x})=\pi^0_1({\bf x})$ and $\pi_2({\bf x})=\pi^0_2({\bf x})$.
	We remark that by a delta method calculation, the asymptotic variance of $\widehat{\theta}=\frac{\widehat{\alpha}_1-\widehat{\alpha}_2}{\widehat{\beta}_2-\widehat{\beta}_1}$ is   $\frac{V(\theta)}{(\beta_2-\beta_1)^2}$, which is proportional to the limit of $n R(\pi_1,\pi_2)$. Thus,  minimizing the latter limit is equivalent to  minimizing the asymptotic variance of  the estimator of the intersection point $\theta$.

	We return to  general dimension $p$.	Set ${\bs \beta}_k^t:=(\beta_{k,1},\ldots,\beta_{k,p})$, $k=1,2$,  ${\bs \beta}_{k,-1}^t:=(\beta_{k,2},\ldots,\beta_{k,p})$ and 
	without loss of generality {assume that $X_1$ is a continuous variable } {given the rest of the covariates}, and that   ${\beta}_{2,1}>{\beta}_{1,1}$. For given ${\bf x}_{-1}=(x_2,\ldots,x_p)$ define $\theta_1:=\theta_1({\bf x}_{-1}):=\frac{\alpha_1-\alpha_2+({\bs \beta}_{1,-1}-{\bs \beta}_{2,-1})^t {\bf x}_{-1} }{\beta_{2,1}-\beta_{1,1}}$. If $\theta_1\in[0,1]$, treatment  $T_1$ is better for a covariate vector $\bf x$ satisfying $x_1 < \theta_1=\theta_1({\bf x}_{-1})$, and otherwise $T_2$ is better. The limit of the regret is given in the following theorem.
	
	\begin{theorem}\label{thm:2}
		Assume model \eqref{eq:model_two} and  $\nu_1,\nu_2 >0$, then
		\begin{equation}\label{eq:asymp_reg_two}
			\lim_{n \to \infty} n R(\pi_1,\pi_2)=\frac{1}{2(\beta_{2,1}-\beta_{1,1})}\int_{[0,1]^{p-1}} V(\theta_1,{\bf x}_{-1}) f(\theta_1,{\bf x}_{-1})I\{\theta_1 \in [0,1]\}d{\bf x}_{-1}.
		\end{equation}
	\end{theorem}
	Thus, asymptotically, $n R(\pi_1,\pi_2)$ is proportional to the integral of the variance of the estimate of the intersection curve $\theta_1({\bf x}_{-1})$, with  weights proportional to the density at the intersection points, and    by \eqref{eq:asymp_reg_two}   $\lim_{n \to \infty} n R(\pi_1,\pi_2)$ is minimized by the allocation $\pi^0_1({\bf x}), \pi^0_2({\bf x})$     since this allocation minimizes $V$ uniformly.
	By Theorems \ref{thm:mult_K_2} and \ref{thm:one_dim_lim} we thus have 
	\begin{theorem}
		\label{th:optp}
		The allocation $\pi^0_1({\bf x}), \pi^0_2({\bf x})$ is  optimal in the sense of  minimizing $\lim_{n \to \infty} n R(\pi_1,\pi_2)$ for any number of covariates $p$. 
		
		Moreover,
		for any $\ve>0$ there exists $C>0$ such that for any allocation functions $\pi_1,\pi_2$ with $\nu_1,\nu_2>0$
		\[
		R(\pi_1,\pi_2) \ge R(\pi_1^0,\pi_2^0) - \frac{C}{n^{3/2-\ve}}.
		\] In particular, the optimal design that minimizes $R(\pi_1,\pi_2)$  can improve the allocation  $\pi_1^0,\pi_2^0$ by at most  ${C}/{n^{3/2-\ve}}$.
	\end{theorem}
	The error term is meaningful as it is of smaller order than the regret, which  is $1/n$ according to Theorem  \ref{thm:2}. 
	
	{In model \eqref{eq:model_two} separate regressions models are considered for each of the treatments and we used the optimal (best linear unbiased estimators) estimators for each model, which are the least squares estimators, to compute $\hat{\delta}$. Alternatively, one can consider a joint regression model 
		$Y=\sum_{k=1}^2  \{\alpha_k+{\bs \beta}_k^t{\bf X}\}{\mathds{1}_k}{+\ve}$, where  ${\mathds{1}_k}$ is the indicator of treatment $k$. Since we did not assume a common variance among the treatments, the variance of $\ve$ in this model depends on ${\mathds{1}_k}$. If the variances $\sigma_1^2,\sigma_2^2$ are known or can be estimated accurately enough, the weighted-least squares estimator with weights proportional to $1/\sigma_1^2,1/\sigma_2^2$ is optimal. This amounts to making the variances equal, and then the above results imply that a balanced design (i.e., $\pi_1(\x)=1/2$) is asymptotically optimal. }
	
	\subsection{A numerical example\label{sec:2_numeric_example}}
\begin{figure}[H]
		\subfigure[The regression lines]{%
			\includegraphics[width=0.45\textwidth]{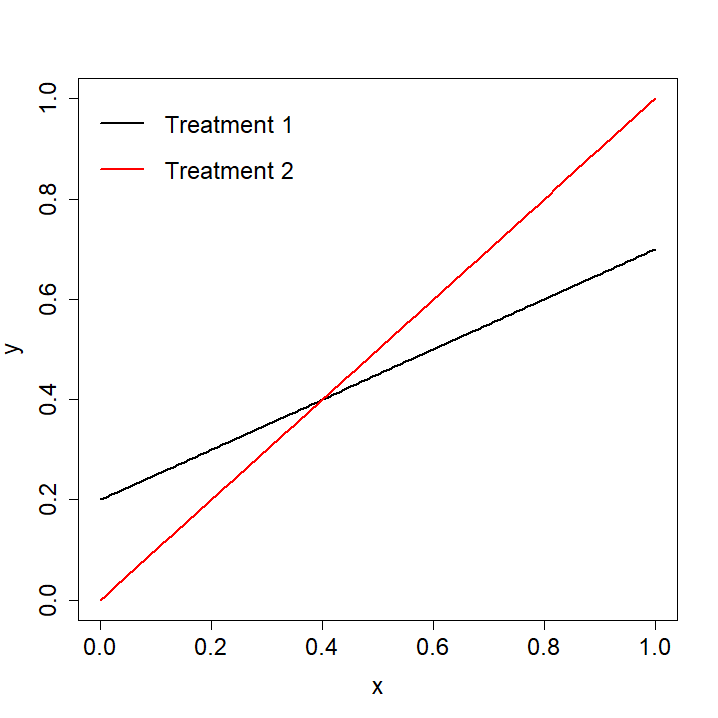}
		}
		\subfigure [Actual and ideal regret  - normal residual]{%
			\includegraphics[width=0.45\textwidth]{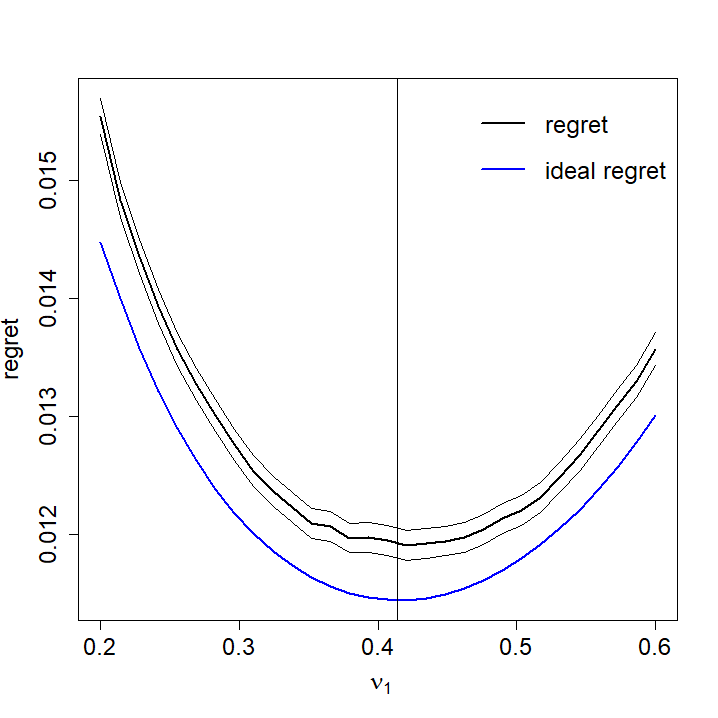}
		}\\
		\subfigure[Actual and ideal regret - exponential residual]{%
			\includegraphics[width=0.45\textwidth]{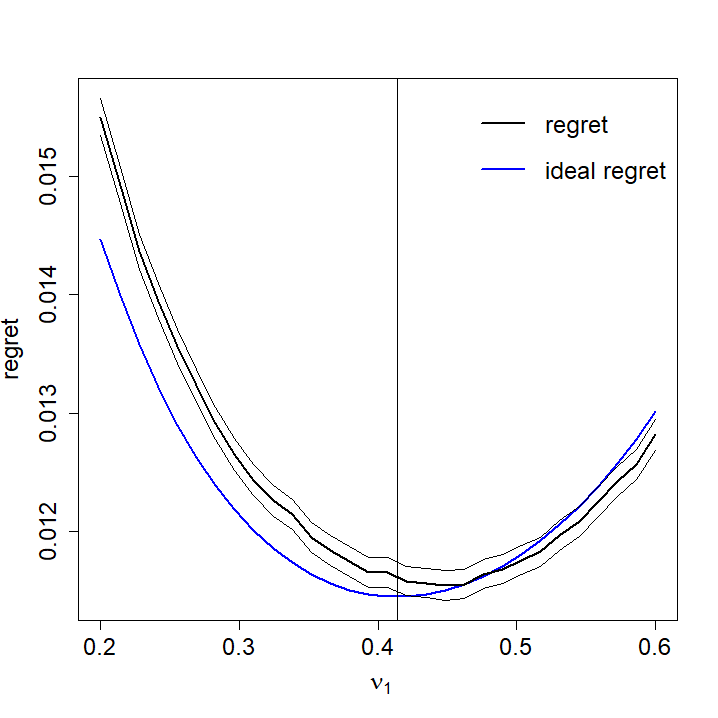}
		}
		\subfigure[Ideal regret and asymptotic regret]{%
			\includegraphics[width=0.45\textwidth]{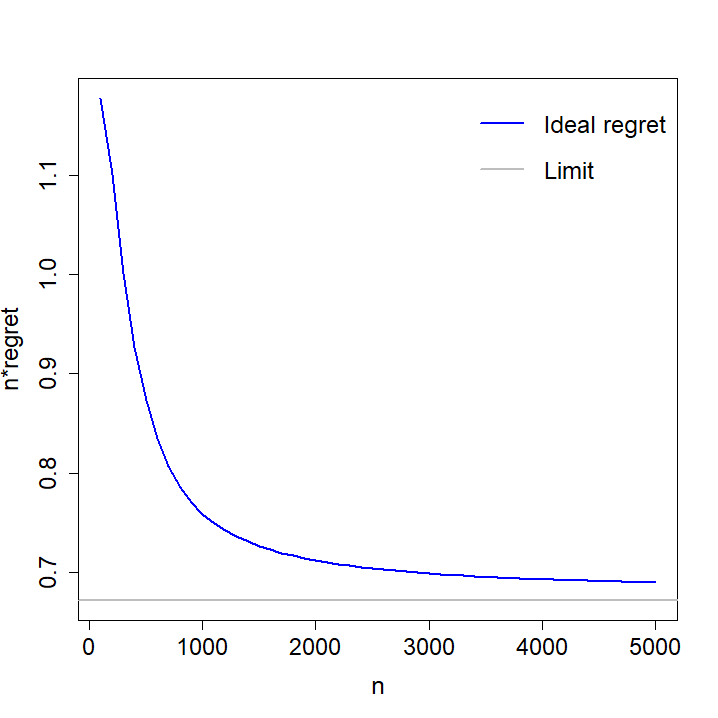}
		}	
		\caption{\footnotesize Figure (a): the regression lines. Figures (b) and (c) show the regret (computed by a simulation with 95\% confidence intervals), and the ideal regret for $n=100$ and different allocation ratios, where in (b) the residual is normal and in (c) it is exponential (centered). The regret is calculated for the allocation function $\pi_1(x)=\nu_1$ where $\nu_1$ varies from 0.2 to 0.6. The optimal allocation is marked by the vertical line. Figure (d) shows $n$ times  the ideal regret (blue line) and its limit (gray line).}\label{fig:regression_model}
	\end{figure}

	To illustrate the asymptotic approximations of the regret and 
	the optimal allocations derived  in Section \ref{sec:approx_regret} 
	consider an example of model \eqref{eq:model_two} with  $p=1$, 
	$\alpha_1=0.2,\, \alpha_2=0,\, \beta_1=0.5,\,\beta_2=1,\,\sigma_1^2=0.1,\,\sigma_2^2=0.2$,
	and $X\sim U[0,1]$. The optimal allocation is $\pi_1({x}) = \nu_1=0.414$. {Notice that this allocation rule does not depend on $x$.} The $\sigma$'s are chosen such that $R^2$ in both regressions is about $0.2$. Figure \ref{fig:regression_model}(a) shows the regression lines. For $x$ smaller than $\theta=0.4$, treatment 1 is better and otherwise treatment 2 is preferred. The regret \eqref{eq:reg2} is  evaluated using simulations (with $10^5$ replications) and is compared to the ideal regret \eqref{eq:reg_approx}. Two scenarios are considered for the residuals in the regression model: normal and  centered exponential. Figures (b) and (c) show the regret and ideal regret for $n=100$ under the allocation function $\pi_1(x)=\nu_1$ {(constant in $x$)} where $\nu_1$ varies from 0.2 to 0.6. The allocation optimizing the ideal regret is marked with a vertical line. While there is a slight deviation between the ideal and actual regret, it seems that the optimal allocations, with respect to both, are very close.   
	Thus, the approximation of Theorem \ref{thm:mult_K_2} works well also for small $n$ in this example. { Figures (b) and (c) also demonstrate the relative robustness of the optimal deign. That is, if instead of implementing the optimal design $\pi_1({x}) = 0.414$, one uses $\pi_1({\bf x}) = \nu_1$ where $\nu_1$ varies from 0.2 to 0.6, the regret loss is relatively small. For example, with $\pi_1(x) = \nu_1$ where $\nu_1=0.3$ and $0.5$, the ideal regret  is larger than the optimal regret by 6\% and 2.9\%, respectively.}      
	Figure \ref{fig:regression_model} (d) shows 	$nR_I(\pi_1,\pi_2)$ of the balanced design ($\pi_1=\pi_2=1/2$) for large $n$ and its limit, see Theorem \ref{thm:one_dim_lim}. It is shown that the regret slowly converges to its limit, and it is close to the limit only when $n$ is quite large.

	\subsection{Polynomial Regression\label{sec:polynomial}}
	We now consider the model of \eqref{eq:model_two} where the functions $g_k(x)$ are assumed to be  polynomials of a single continuous covariate $X$ with density $f$, i.e.,
	$g_k(X)=\alpha_k +\sum_{j=1}^J\beta_{jk} X^j$,
	where $J$ denotes the degree of the polynomial. {We will show that $\pi_1^0(x),\pi_2^0(x)$ are  asymptotically optimal also in this setting.}
	
	Let $\theta_1<\cdots<\theta_L$, where $L \le J$, be the crossing points of $g_1(x)$ and $g_2(x)$, and define $\theta_0=0$ and $\theta_{L+1}=1$. Assume without loss of generality that $g_1(x)>g_2(x)$  for $x\in (\theta_\ell,\theta_{\ell+1})$  if $\ell$ is odd, and the reverse inequality holds when $\ell$ is even. Therefore, the regret is
	\begin{multline*}
		R(\pi_1,\pi_2)=\sum_{\ell\text{ is odd}}\int_{\theta_\ell}^{\theta_{\ell+1}} P(\widehat{g}_2(x) > \widehat{g}_1(x))(g_1(x) - g_2(x))f(x) dx\\+\sum_{\ell\text{ is even}}\int_{\theta_\ell}^{\theta_{\ell+1}} P(\widehat{g}_1(x) > \widehat{g}_2(x))(g_2(x) - g_1(x))f(x) dx.
	\end{multline*}
	
	Notice that Theorems \ref{thm:mult_K_2} and \ref{thm:2}, which provide results for multivariate ${\bf X}$, do not apply in this case because the random vector ${\bf X}=(X,X^2,\ldots,X^J)$ does not have a joint density.  However, using the arguments of Theorem \ref{thm:mult_K_2} (under similar notation and regularity conditions) it can be shown that $n^{3/2-\ve} |R(\pi_1,\pi_2)-R_I(\pi_1,\pi_2)|\to 0$, for any $\ve>0$, where
	\begin{equation}\label{eq:polya}
		R_I(\pi_1,\pi_2)=\int_0^1 \Phi\left(\frac{-\sqrt{n}|g_2(x) - g_1(x)|}{\sqrt{V(x)}}\right)|g_2(x) - g_1(x)| f(x)dx,
	\end{equation}
	and now ${\bf x}^t=(1,x,\ldots,x^J)$, ${\bs Q}_k = \int_{\mathbb{R}^p} {\bf x} {\bf x}^t f_k({x}) d{x}$ and ${\bs \Sigma}_{k}
	=\frac{1}{\nu_k}\sigma^2_k {\bs Q}_k^{-1}
	, ~~k=1, 2$,
	and  $V(x)= {\bf x}^t\left( {\bs \Sigma}_{1}+ {\bs \Sigma}_{2} \right){\bf x}$, which is the asymptotic variance of $\widehat{g}_1(x)-\widehat{g}_2(x)$.

	Theorem \ref{th:optp} implies that the design where $\pi_1(x)=\pi_1^0(x)$ minimizes  $V(x)$ uniformly over all $x$'s. It follows that this design is asymptotically optimal also for this problem.

	By Proposition \ref{prop:poly} (in the appendix) we have 
	\[
	\lim_{n \to \infty} n R(\pi_1,\pi_2)=\sum_{\ell =1}^L \frac{f(\theta_\ell)V(\theta_{\ell})}{2|({\bs \beta}_2-{\bs \beta}_1)^t{\bs \zeta}_\ell|},
	\]
	where ${\bs \zeta}_\ell:=(1,2\theta_\ell,...,J\theta_\ell^{J-1})^t$.
	
	A careful inspection of the proofs of the above results shows that  the optimality of the design with $\pi_1(x)=\pi^0_1(x)$  generalizes   to regression functions of the form $g_k(X)=\sum_j \beta_{jk} h_{j}(X)$ for any  functions $h_{j}$ having bounded derivatives.

	\section{$K$ Treatments and one Covariate}\label{sec:K_treatments}
	\subsection{Regret and ideal regret}\label{sec:Knorm}
	
	We consider the case $p=1$ and 
	$g_k({  X})=\alpha_k+{  \beta}_k  {  X}, ~Var(Y|{  X},T_k)=\sigma_k^2,~k=1,\ldots,K  $.
	As in Section \ref{sec:two_treatments},  we assume the existence of the moment generating function of $\ve$ when conditioned on $({X},T_k)$. We also assume that ${X}$ is continuous with density $f({x})$ and is supported on $[0,1]$. Let $\theta_k$, $k=0,\ldots,K$ be in increasing order, and such that treatment $k$ is  best in the interval $(\theta_{k-1}, \theta_k)$; we assume that each treatment is best in  some open interval, or equivalently that the intervals are nonempty.
	
	Equations \eqref{eq:regretdirect} and \eqref{eq:IR} for the regret and ideal regret imply immediately
	\[
	R(\pi_1,\ldots,\pi_K)=\sum_{k=1}^K\sum_{m=1}^K\int_{\theta_{m-1}}^{\theta_{m}}P\left(\widehat{g}_k(x) > \max_{\ell \ne k} \widehat{g}_\ell(x) \right)\left[ g_{m}({x}) - g_k({x})\right]f({x})d{x}
	\]
	and the ideal regret is
	\begin{align} \label{eq:regretdirectdetailmdy}
		&R_I(\pi_1,\ldots,\pi_K)\nonumber\\&=\sum_{k=1}^K\sum_{m=1}^K\int_{\theta_{m-1}}^{\theta_{m}}\left(\int \prod_{\ell=1,\ldots,K,\ell\not=k} \Phi\left(\frac{z\xi_k(x)+\sqrt{n}[g_k({x})-g_\ell({x})]}{\xi_\ell(x)}\right)\varphi(z)dz\right)\left[ g_{m}({x}) - g_k({x})\right]f({x})d{x}.
	\end{align}

	Similar to the results in Section \ref{sec:two_treatments}, we give asymptotic results on the rate of convergence of the regret to the ideal regret $R_I$ as well as the limit of the regret.
	Given  allocation functions $\pi_k(x)$, we have by \eqref{eq:def_Q} that the distribution of the estimated regression coefficients $(\widehat \alpha_k,\widehat \beta_k)$ is approximately bivariate normal with means $(\alpha_k,\beta_k)$ and 
	covariance
	\begin{equation}\label{eq:Sigma_beta}
		{\bs \Sigma}_{k}=\frac{\sigma^2_k}{\nu_k}\left(
		\begin{array}{ll} \frac{\tau^2_k+\mu_k^2}{\tau^2_k}&
			\frac{-\mu_k}{\tau^2_k}\\
			\frac{-\mu_k}{\tau^2_k}&\frac{1}{\tau^2_k}
		\end{array}\right)\,,
	\end{equation}
	where we assume positivity of $\nu_k:=\int f({x}) \pi_k({ x}) d{ x}$,  and denote the mean and variance of the covariate in group $k$ by $\mu_{k}:=\int_{0}^{1}x f_k(x)\,dx$ and $\tau_k^2:=\int_{0}^{1}(x-\mu_k)^2 f_k(x)\,dx$, respectively,   $k=1,\ldots,K$. Parallel to Theorems \ref{thm:mult_K_2} and \ref{thm:one_dim_lim} we have
	\begin{theorem}
		\label{thm:K_one_dim1}
		Under the assumptions in the beginning of Section \ref{sec:Knorm} we have for any $\ve>0$
		\[
		\lim_{n \to \infty} n^{3/2-\ve} \left|R(\pi_1,\ldots,\pi_K)-R_I(\pi_1,\ldots,\pi_K)\right| = 0.
		\]
	\end{theorem}
	\begin{theorem}
		\label{thm:K_one_dim_lim}
		Under the assumptions in the beginning of Section \ref{sec:Knorm}, we have 
		\begin{equation}
			n\lim_{n \to \infty}R(\pi_1,\ldots,\pi_K)=\sum_{m=1}^{K-1}\frac{V_m(\theta_m)f(\theta_m)}{2|\beta_{m+1}-\beta_{m}|}, \label{eq:asymptotic_regret_k_trt}    
		\end{equation}
		where  $V_m(x)=(1,x)\left( {\bs \Sigma}_{m}+ {\bs \Sigma}_{m+1} \right)\binom{1}{x}$.
	\end{theorem}
	{Notice that the assumption that the intersection points $\theta_1,\ldots,\theta_K$ are in $(0,1)$ implies that $\beta_{m+1}\ne\beta_{m}$ for $m=1,\ldots,K-1$.} Using \eqref{eq:Sigma_beta} it is possible to write $V_m(\theta_m)$ explicitly as follows
	\[
	V_m(\theta_m)=\frac{\sigma_{m}^2}{\nu_{m}}\left[1+\frac{(\theta_m-\mu_{m})^2}{\tau_{m}^2}\right]+\frac{\sigma_{m+1}^2}{\nu_{m+1}}\left[1+\frac{(\theta_m-\mu_{m+1})^2}{\tau_{m+1}^2}\right].
	\]
	Theorem \ref{thm:K_one_dim_lim} implies that asymptotically the optimal allocation problem reduces {to minimize a weighted average of the variances of  $\widehat{\theta}_1,\ldots,\widehat{\theta}_{K-1}$ (see the discussion after Theorem \ref{thm:one_dim_lim}), which are the estimates of the intersection points. The weights depend on the $\beta$'s and on the density $f$ at the intersection points.} 
	Note also that $V_m(\theta_m)$ is the sum of the asymptotic variances of $\widehat \alpha_m + \widehat \beta_m \theta_m$ and  $\widehat \alpha_{m+1} + \widehat \beta_{m+1} \theta_m$. 	
	Theorem \ref{th:optp} continues to hold for $K$ treatments and thus, optimizing the ideal regret approximately optimizes the regret itself.
	\subsection{Approximate numeric optimization of allocation rules}
	\label{sec:numopt}
	In this section we demonstrate the utility of (\ref{eq:regretdirectdetailmdy}) in computing optimal designs when $K>2$, and the reduction in the regret that can be achieved. {The focus here is on the ideal regret. By Theorem \ref{thm:K_one_dim1} and the numerical results of Section \ref{sec:2_numeric_example},  a reduction in the ideal regret implies a similar reduction in the regret itself.}  
	
	To numerically find allocation rules minimizing the regret (\ref{eq:regretdirectdetailmdy}), we can use a finite dimensional parametrization of  the set of allocation probability functions $\pi_k(x), \,k=1,\ldots,K$.  One option is to use suitably rescaled polynomials of degree  $M\geq 0$: Let ${\bf A}=(a_{k,m})$ denote a $(K-1) \times (M+1)$ matrix and set 
	$h_k(x)=\exp(\sum_{m=0}^M a_{k,m} x^m)$ for 
	$k=1,\ldots,K-1$.  Then, the allocation probabilities are  defined by 
	$\pi_k(x)=\frac{h_k(x)}{1+\sum_{j=1}^{K-1} h_j(x)},\, k=1,\ldots,K-1$ and $\pi_K(x)=1-\sum_{j=1}^{K-1}\pi_j(x)$.
	Note that  $\mathbf{A}=\mathbf{0}$ corresponds to equal allocation $\pi_k(x)=1/K,k=1,\ldots,K$. Now, we can  optimize the ideal regret by plugging $\pi_k(x)$ into \eqref{eq:regretdirectdetailmdy} and numerically minimizing with respect to ${\bf A}$ (e.g., by using  the algorithms provided with the R function optim).  
	
	For $M=0$ we obtain (approximately) \textit{optimal fixed allocation} probabilities that do not depend on $x$.  The larger $M$ the more flexible allocation functions in $x$ are fitted. However, because the target function (\ref{eq:regretdirectdetailmdy}) is not necessarily convex, we cannot guarantee that the numerical optimization will converge to the global minimum regret. Because the spaces of allocation probability functions are nested for increasing $M$, a heuristic optimization strategy is to iteratively  optimize (\ref{eq:regretdirectdetailmdy}) starting with $M=0$ and then subsequently increase $M$  using the solution of the last iteration step as starting values for the next iteration (adding a column of  zeroes to $\mathbf{A}$).  Note that this approach can also be used to optimize allocation rules minimizing the asymptotic regret (\ref{eq:asymptotic_regret_k_trt}) instead of the ideal regret (\ref{eq:regretdirectdetailmdy}).
	
	To illustrate this approach we consider a specific setting where $n=200$, $K=3$, $(\alpha_1, \alpha_2, \alpha_3)=( 0.0, -0.1, -1.2)$ and $(\beta_1,\beta_2,\beta_3)=(0.2, 0.5, 2.0)$; the left panel of Figure \ref{fig:exampleuniform} plots the regression lines. We looked at two options for the variances: heteroscedastic - where $(\sigma_1,\sigma_2,\sigma_3)=\frac{1}{\sqrt{12}}(\beta_1,\beta_2,\beta_3)$ (so that $\sigma^2_k=Var(\beta_k X)$ when $X$ is uniform for each $k$, i.e., $R^2\approx 1/2$ for each regression) and homoscedastic - where $\sigma_1=\sigma_2=\sigma_3=\frac{\beta_2}{\sqrt{12}}$. We consider four settings for the distribution of $X$: Uniform(0,1) = Beta(1,1), Beta(0.5,0.5) (U-shaped), Beta(2,2) (symmetric, unimodal), Beta(2,5) (asymmetric). 
	
	{Starting with the balanced design, Table \ref{tab:reg} shows the percentage of regret reduced by using the optimized fixed allocation design ($M=0$) and by the optimized design with $M=4$. 
		The allocation probabilities when $M=0$ are also given. It is clear from the table that the regret depends heavily on the distribution of $X$ and on the $\sigma$'s. 
		The table demonstrates that the reduction of the regret is quite significant. Furthermore, the optimal design where the allocation depends on $x$ improves significantly over the fixed design ($M=0$). An exception is the homoscedastic case and $X\sim Beta(2,5)$, in which the density function is small when $g_3(x)$ is maximal (see the left panel of Figure \ref{fig:exampleuniform}) and additionally $\sigma_3$ is small compared to $\beta_3$. This makes this case close to $K=2$ where fixed allocation designs are optimal.    
		Another finding, which is not given in the table, is that the difference, in terms of the regret, between $M=1$ and $M=4$ is small; it is at most 2.3\% in the eight scenarios we considered.    }

	\begin{table}[ht]
		\centering
		\begin{tabular}{|c|c|c|c|c|c|c|c|c|c|c|}\hline
			& \multicolumn{5}{c}{heteroscedastic}&\multicolumn{5}{|c|}{homoscedastic} \\\hline 
			& \multicolumn{2}{c}{\% reduction}&\multicolumn{3}{|c|}{fixed alloc. ($M=0$)} & \multicolumn{2}{c}{\% reduction}&\multicolumn{3}{|c|}{fixed alloc. ($M=0$)} \\\hline 
			distribution & $M=0$ & $M$=4 &$\pi_1$&$\pi_2$&$\pi_3$& $M=0$ & $M$=4 &$\pi_1$&$\pi_2$&$\pi_3$\\ \hline
			U(0,1) &24.0 &  40.9 &0.12& 0.33 & 0.55 & 11.7 & 25.4  & 0.40 & 0.43 & 0.17 \\
			Beta(0.5,0.5)& 25.1 & 41.9& 0.13 & 0.34 & 0.53 & 16.5 & 29.6 & 0.41 & 0.44 & 0.15\\
			Beta(2,2)& 20.5 & 34.2 & 0.11 & 0.32 & 0.57 & 7.0 & 18.8 & 0.38 & 0.43 & 0.19 \\
			Beta(2,5)& 11.5 & 22.8 & 0.18 & 0.45 & 0.37 & 21.3 & 25.8 & 0.45 & 0.46 & 0.09 \\ 
			\hline 
		\end{tabular}	
		{\caption{Results for the setting mentioned in the text under different distributions of $X$ and for the heteroscedastic and homoscedastic models: the percentage of reduction achieved by the optimized designs with $M=0$ and $M=4$ relative to a balanced design, and the fixed allocation probabilities for $M=0$.\label{tab:reg}}}
	\end{table}

	\subsection{A lower bound for the asymptotic regret  (\ref{eq:asymptotic_regret_k_trt}) }\label{sec:lowerbound}
	By  (\ref{eq:regretdirectdetailmdy}), \eqref{eq:Sigma_beta}, and \eqref{eq:asymptotic_regret_k_trt} (see also \eqref{eq:def_Q}) both the ideal and the asymptotic regret depend on the functions  $\pi_k(x)$ only via the quantities $\nu_k$, $\mu_k, \tau_k^2,k=1,\ldots,K$. Therefore, we can minimize (\ref{eq:regretdirectdetailmdy}) and (\ref{eq:asymptotic_regret_k_trt}) in these parameters to obtain a lower bound for the regret. 
	Let $\mu:=\int_{0}^{1} xf(x) \,dx$,  $\tau^2:=\int_{0}^{1} (x-\mu)^2f(x)\,dx$ and recall that  
	$f(x)=\sum_{k=1}^K\nu_k f_k(x)$. 
	Then, by the relations for  central moments of  mixture distributions we have
	\begin{equation} \label{eq:const}
		\sum_{k=1}^K\nu_k=1,\quad \sum_{k=1}^K\nu_k \mu_k=\mu,\quad \sum_{k=1}^K\nu_k (\tau_k^2+\mu_k^2)=\tau^2+\mu^2.  
	\end{equation}
	Minimizing (\ref{eq:regretdirectdetailmdy}) or (\ref{eq:asymptotic_regret_k_trt}) in  $\nu_k$, $\mu_k, \tau_k,k=1,\ldots,K$  subject to the constraints (\ref{eq:const}) and $\nu_k,\tau_k\geq 0,k=1,\ldots,K$ yields a lower bound for the achievable regret. 
	To obtain the minimum regret (instead of a lower bound) one  needs to add additional constraints on  $\mu_k, \tau_k,k=1,\ldots,K$ to restrict optimization to values for which there exist mixture components $f_k$  with weighted sum $f$ that assume these moments.
	
	We now consider an example with a uniform $f$ and $K=3$, and minimize the asymptotic regret (\ref{eq:asymptotic_regret_k_trt}). Since  now $f(x) \le 1$ the constraint  $f(x)=\sum_{k=1}^3 \nu_k f_k(x)$  implies that  $\nu_k f_k(x)\leq 1,\,  k=1,2,3,\, x\in[0,1]$. By Lemma 1 below this implies that the variance of  $f_k(x),k=1,2,3$ is bounded from below by $\nu_k^2/12$. 
	
\begin{figure}
		\begin{tabular}{cc}\includegraphics[width=7cm]{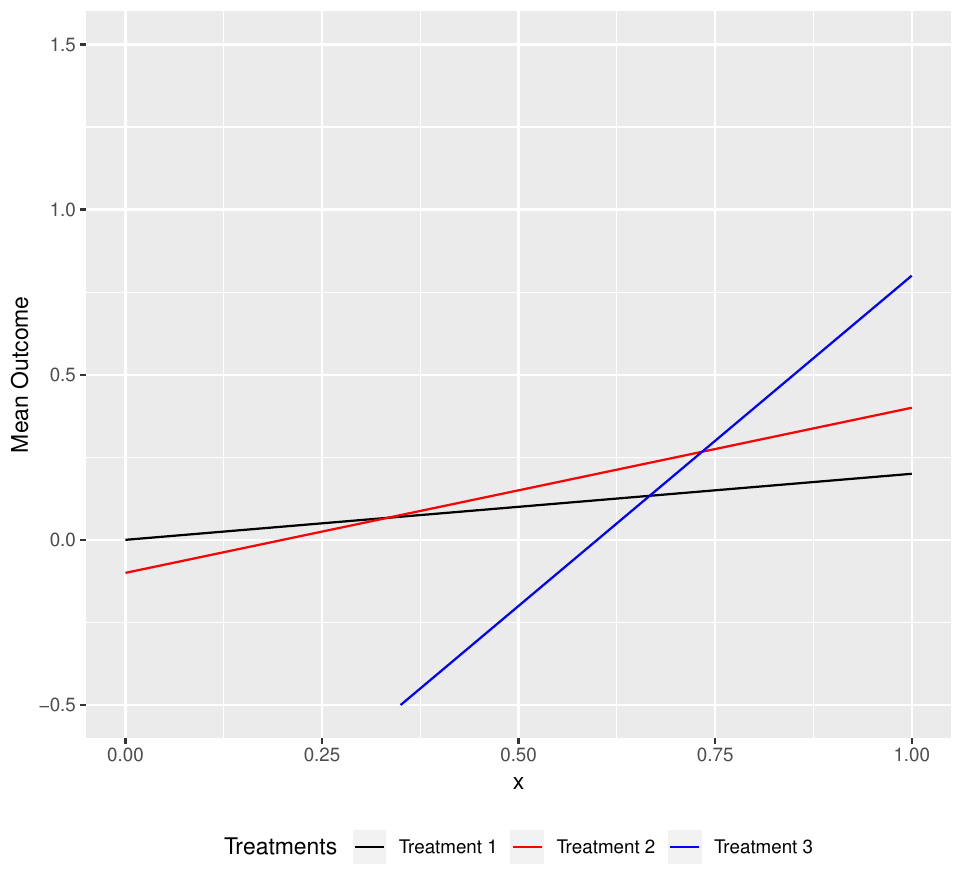}&
			\includegraphics[width=7cm]{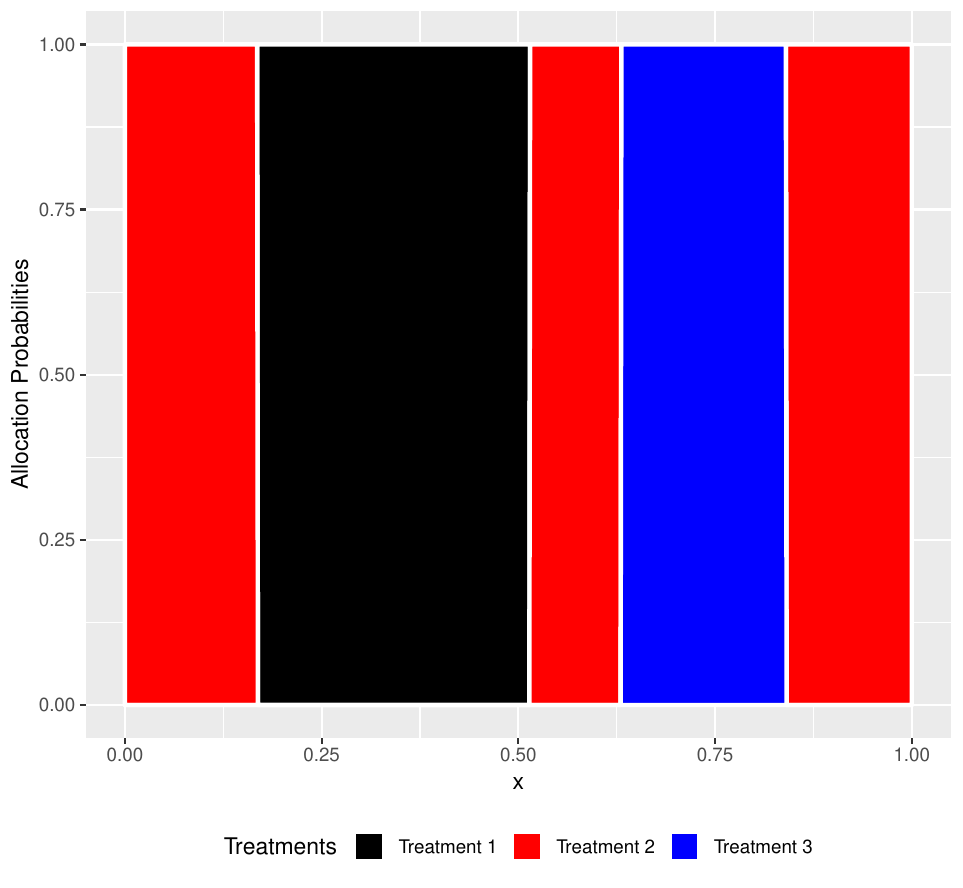}
		\end{tabular}
		\caption{The functions $g_k$ for the three treatments (left) and the asymptotically optimal allocation probabilities minimizing \eqref{eq:asymptotic_regret_k_trt} (as a stacked area chart, right) for the example in Section \ref{sec:lowerbound}. For instance, for $x$'s under the black area, treatment 1 is assigned with probability $\pi_1(x)=1$. }
		\label{fig:exampleuniform}
	\end{figure}

	
	As an example consider the setting of Section \ref{sec:numopt}
	where $(\alpha_1, \alpha_2, \alpha_3)=( 0.0, -0.1, -1.2),\,  (\beta_1,\beta_2,\beta_3)=(0.2, 0.5, 2.0)$, and the residual variances $\sigma_k^2,\, k=1,2,3$ are equal across treatment groups.
	Figure \ref{fig:exampleuniform} shows the
	resulting scenario and optimal allocation probabilities (whose computation is explained below). Then, numerically (see details below) minimizing \eqref{eq:asymptotic_regret_k_trt} in  $\nu_k$, $\mu_k, \tau_k,k=1,\ldots,K$  subject to the constraints (\ref{eq:const})  as well as 
	$\tau_k^2\geq\nu_k^2/12,\, k=1,2,3$
	gives $\nu_k=0.346,\, 0.444, \,0.210, \,
	\mu_k=0.342,\, 0.512,\, 0.735, \,
	\tau^2_k=0.00997,\, 0.132,\, 0.00369 $,  respectively,
	and the minimized  value of  \eqref{eq:asymptotic_regret_k_trt} is $12.128$. {At this point we know it is a lower bound to \eqref{eq:asymptotic_regret_k_trt} since only part of the constraints were applied.  However, we now demonstrate a design which achieves this lower bound. }
	
	We note that $\tau_k^2=\nu_k^2/12$,  for $k=1,3$. Applying Lemma 1 for $c=1/\nu_k$  for $k=1,3$, the only density $f_k(x)$ with mean $\mu_k$ and variance $\tau_k^2$ and satisfying the constraint $f_k(x)\leq 1/\nu_k,\,x\in[0,1]$ is the uniform distributions on $[\mu_k-\nu_k/2 ,\, \mu_k+\nu_k/2]$ with densities $f_k(x)=1_{[\mu_k-\nu_k/2 ,\,\, \mu_k+\nu_k/2]}(x)/\nu_k$.
	Furthermore, we set
	$$f_2(x)=\frac{1-\nu_1 f_1(x)-\nu_3 f_3(x)}{\nu_2}\,.$$
	Because the supports of $f_1,f_3$ are disjoint and $f_k(x)\leq 1/\nu_k,\,k=1,3$, $f_2(x)$ is a valid density. 
	By construction, the mean and variance of $f_2(x)$ are $\mu_2,\tau^2_2$\, and the three densities satisfy $f(x)=\sum_{k=1}^3 \nu_k f_k(x)$. 
	The resulting allocation probabilities $\pi_k(x)=f_k(x)/\nu_k$ achieve the lower bound computed above, and hence, assuming the above minimization algorithm yielded the correct minimum (which we verified in several ways, including the fact that it is unique, see below), it is an optimal allocation (see Figure \ref{fig:exampleuniform} for a plot). It follows that in this case the above constraints suffice.  Furthermore, as  $f_1,f_3$ are unique, they are a unique optimal allocation (under the assumption that the minimization in $\nu_k$, $\mu_k, \tau_k,k=1,\ldots,K$ gives a unique solution). It also follows that the optimal solution leads to a deterministic allocation. 
	
	Some comments: (1)  To implement the numerical optimization of (\ref{eq:asymptotic_regret_k_trt}) subject to the constraints, we use the R command optim with the following reparametrization, which allows one to apply optimization algorithms for unrestricted optimization. Let $A_i=(a_{i,1},a_{i,2})\in R^2,\, i=1,2,3$. First we transform the vectors $A_i$ into proportions setting $a'_{ik}=\exp(a_{i,k})/\sum_{m=1}^3\exp(a_{i,m}),\, i,k=1,2,3$ where $a_{i,3} \equiv 0$, and set $\nu_k=a'_{1k},\, \mu_k=a'_{2k} \mu /\nu_k$, and
	$$
	\tau^2_k=\frac{\nu_k^2}{12}+\frac{a'_{3,k}\left(\mu^2+\tau^2-\sum_{m=1}^3(\nu_m\mu_m^2+\nu_m^3/12)\right)}{\nu_k} , \,\,k=1,2,3.
	$$
	As long as the second summand remains positive, this parametrization ensures that $\tau^2_k\geq \nu_k^2/12$ and that the constraints of \eqref{eq:const} hold. To confirm the robustness of the solution, the optimization  {was} performed repeatedly with randomly chosen initial conditions. 
	(2) If the intersection points $\theta_{k,k+1}$ of the  lines $g_k$ are moved closer together by a different specification of the parameters, so will the intervals of $x$'s where all patients are allocated to Treatment 1 or Treatment 3 (the blue and black regions in Figure \ref{fig:exampleuniform}). For scenarios where they overlap, we can no longer use the above approach to obtain an optimal solution and the solutions seem to become more complex. This is the case, for example, for $\alpha_k= 0.0,\, -0.1,\, -.8,\, \beta_k= 0.2,\, 0.5,\, 2.0,\,\,k=1,2,3$.
	
	\paragraph{\bf Lemma 1}  Let $c> 0$. Among all continuous distributions $f$  with (i) mean $\mu$ such that (ii) $f(x)\leq c$ for all $x\in \mathbb{R}$, the uniform distribution on $[\mu-1/(2c), \mu+1/(2c)]$ has the smallest variance. This minimum variance is $1/(c^2 12)$.

	\section{A clinical trial to determine personalized diets}\label{sec:diets}
	
	To demonstrate our approach, we consider a study by \cite{Ebbeling2018} in which three diets were compared in a parallel group design. We fitted a model to the data of that trial and derived an alternative optimal design based on the estimated model parameters.
	
	Subjects were randomly assigned to three diet groups. The diets differ in their carbohydrate content, high (60\%), moderate (40\%) or low (20\%) and were followed for 20 weeks. The total number of subjects allocated to these three treatment groups was 54, 53 and 57, respectively, so the design is (almost) balanced. The main {covariate} is {insulin secretion} (insulin concentration 30 minutes
	after oral glucose).
	The primary outcome in the original trial was averaged total energy expenditure over two measurements, in the middle and the end of the trial.
	Total energy expenditure is  measured by doubly labeled water; for details on this measure see \citet{Hills2014}.
	
	Of the 1685 subjects initially screened, 234  participated in a
	run-in phase 
	which preceded the trial itself. Of these, 164 achieved
	a target of	12\% ($\pm$ 2\%) weight loss and were qualified for the trial; they were randomly
	assigned to one of the three diets.

	To account for the cost of the diets, either in terms of money or  in terms of inconvenience to the dieter, we consider the utility of the diets as an outcome variable which is defined as the energy expenditure minus the diet cost. We assume that   the lower the carbon content of the diet, the higher the cost (see  \citet{Hagberg2019} for a similar definition) but otherwise our choice of the cost for this example is arbitrary.  On the basis of Figure 4 in \cite{Ebbeling2018}, we assume that under treatment $T_k$, $k=1,2,3$,  
	\begin{equation} \label{eq:diets}
		\text{utility}=\alpha_k + \beta_k X + \ve_k - \text{cost}_k,    
	\end{equation}
	where cost$_k$= 0, 150, 300 for $k=1,2,3$, respectively, $X$ is insulin secretion and its distribution is Gamma(3.12, 0.02), and the standard deviations of $\ve_1,\ve_2,\ve_3$ are 190, 150, 130, respectively. The design we propose is based on the assumed values of $\alpha_k$, which are 40, -80, -240, and of $\beta_k$, which are -0.8, 0.1, and 0.8 for $k=1,2,3$, respectively. These values and the parameters of the gamma distribution of $X$ and the standard deviation of $\ve$'s are all based on Figure 4 of \citet{Ebbeling2018}, which plots the energy expenditure as a function of the insulin secretion for the three diets in the data, {and roughly shows the distribution of insulin secretion}. In particular, the parameters of the Gamma distribution are chosen to match empirical quantiles of the $X$'s.
	
	{In order to numerically find the optimal allocation probabilities $\pi_k(x)$ that minimize the ideal regret (\ref{eq:regretdirectdetailmdy}), we use  the approach outlined in subsection \ref{sec:numopt} setting $M=4$ and using the R function optim with the Nelder-Mead method. }

\begin{figure}[ht!]
		\subfigure[the regression models]{%
			\includegraphics[width=0.45\textwidth]{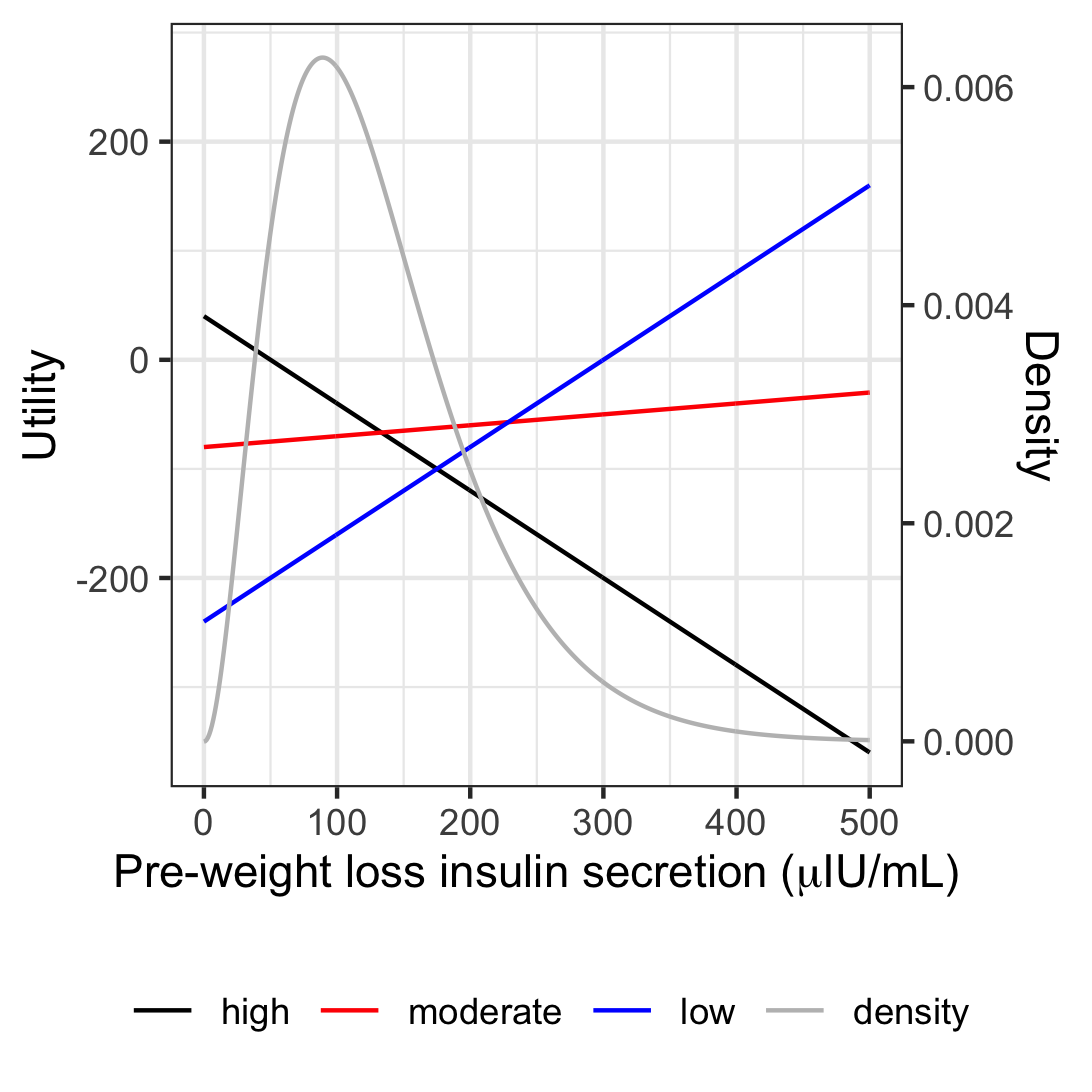}
		}
		\subfigure [optimal allocation for $n=164$]{%
			\includegraphics[width=0.45\textwidth]{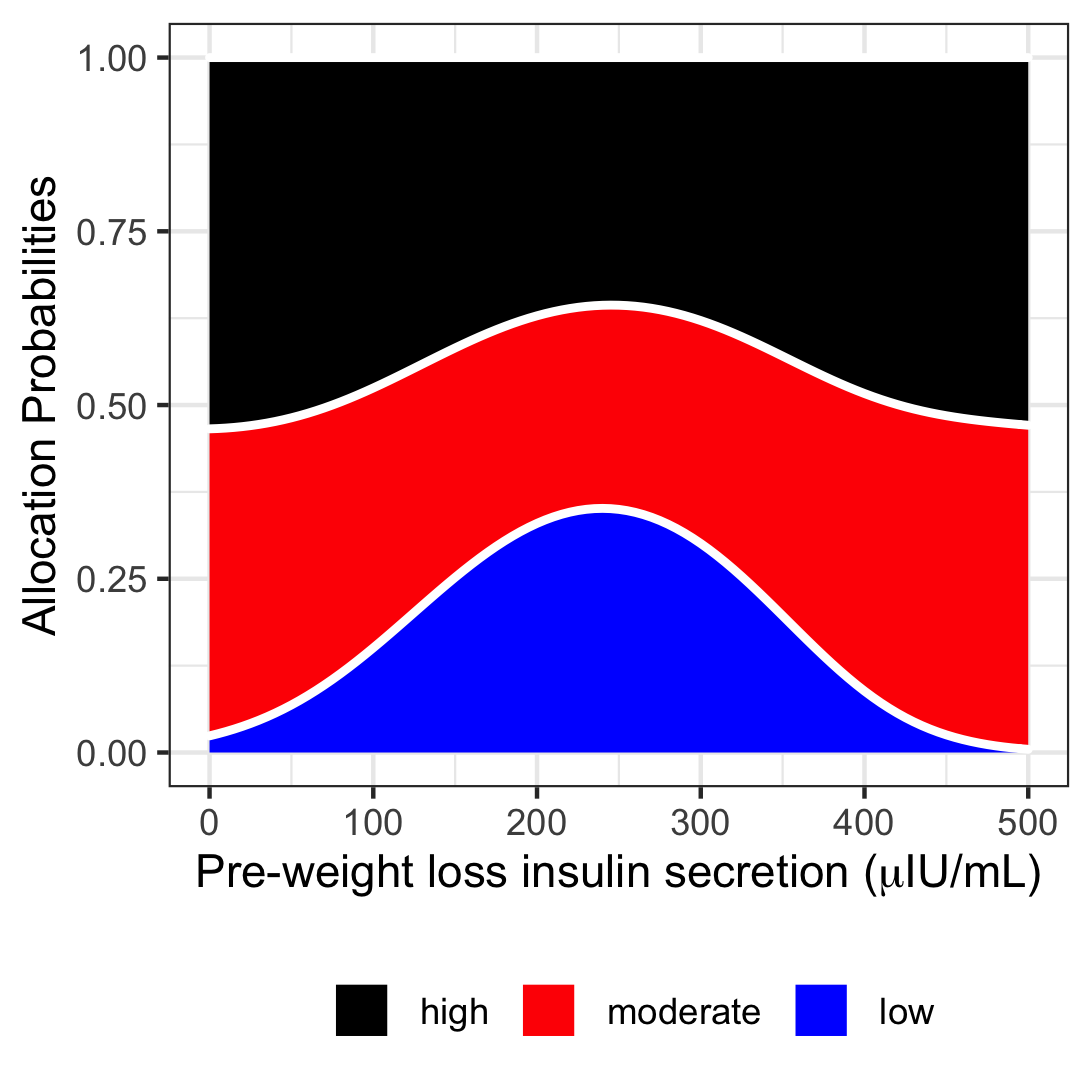}
		}\\
		\subfigure[optimal allocation for $n=1000$]{%
			\includegraphics[width=0.45\textwidth]{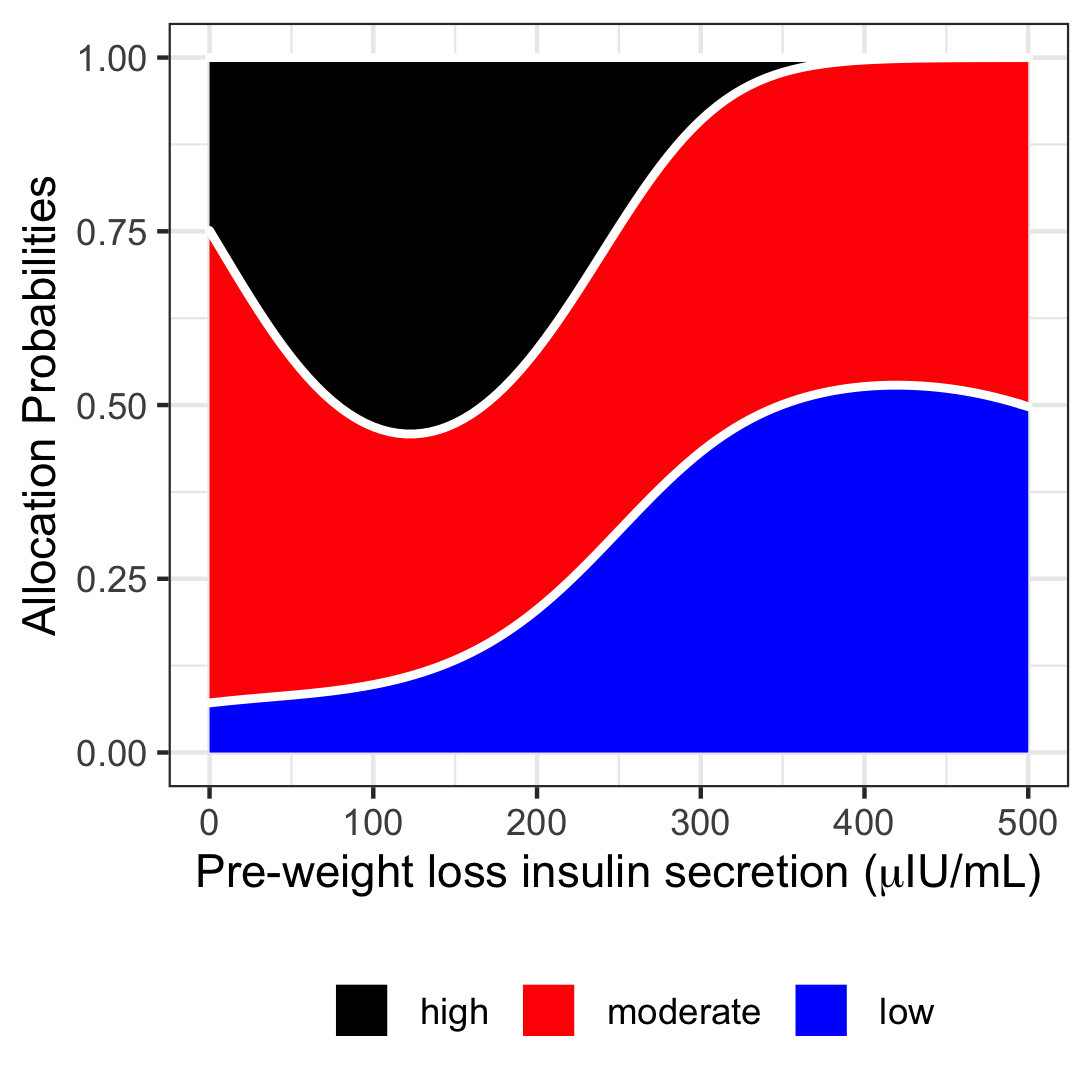}
		}
		\subfigure[optimal allocation for $n\to\infty$]{%
			\includegraphics[width=0.45\textwidth]{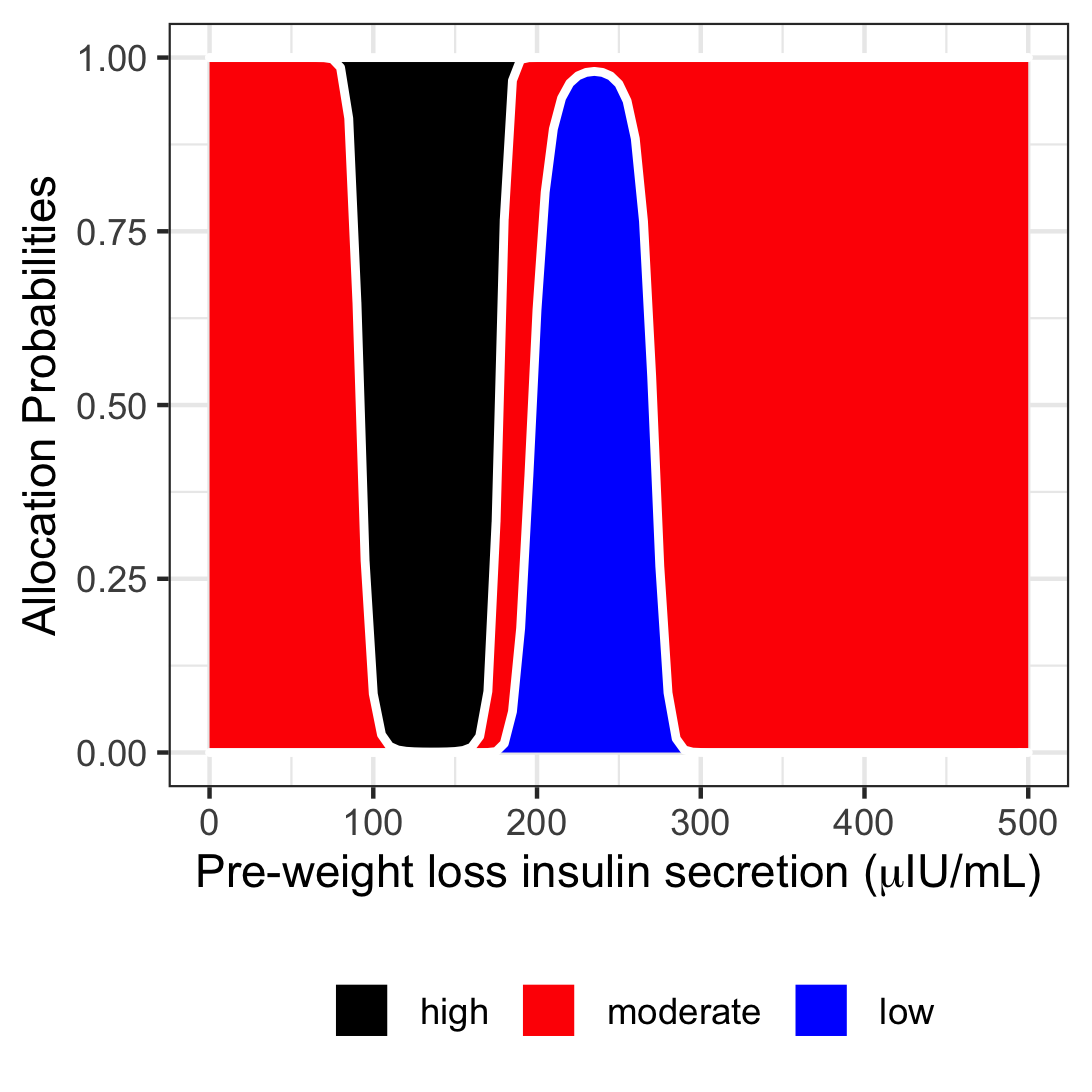}
		}	
		\caption{Figure (a) plots the regression models \eqref{eq:diets} and the assumed density of the covariate. Figures (b)-(d) present the optimal allocation probabilities (as a stacked area chart) when $n=164, 1000,  \infty$.}\label{fig:diets}
	\end{figure}

	Figure \ref{fig:diets}(a) plots the regression lines \eqref{eq:diets} for the three diets and the assumed distribution of the covariate $x$. The high, moderate, and low carbohydrate diets are optimal when $x$ belongs to the intervals $(0, 133),(133, 229),(229, \infty)$, respectively. The difference between the treatments is more pronounced for extreme values of $x$. Hence, people with very small or very large $x$ would benefit more from a personalized treatment choice compared to those with medium values of $x$. 

	Figures \ref{fig:diets}(b)-(d) present the optimal allocation probabilities, as a stacked area chart, when $n=164,\, 1000$, and $\infty$. When the sample size is larger, the optimal allocation rule is closer to being deterministic.  For $n=164$, which was the actual trial size, the optimal trial allocation rule tends to allocate subjects to $T_3$ (respectively, $T_2$)  for large (respectively, small) values of $x$. Compared to equal allocation, which was the actual design used, there is a reduction of about $11\%$ of the regret. A saving of 11\% may be quite significant if the treatments are applied to a large number of future patients. Furthermore, the regret of the balanced design with $n=164$ can be achieved under the optimal design with $n=143$, which represents a reduction of 13\% of the sample size. Taking into account that in this experiment only about 10\% of recruited  patients passed the screening process,  a saving of about 20 subjects amounts to a reduction of about 200 subjects that would have to be recruited and pretested. Also, the duration of 20 weeks, and the cost involved in the experiment, imply that an optimal design is desirable.
	
	To investigate if the reduction in the regret is due to unbalanced allocation or the covariate dependent allocation, we also computed the regret for the optimal fixed allocation that is restricted to allocation probabilities that do not depend on the covariate.  This can be achieved as above by setting $M=0$. We found numerically the optimal fixed design  for $n=164$ leads to reduction of the regret by about 8\%, whereas optimizing with allocation functions that depend on the covariate, leads to a reduction by an additional 3\%. For larger sample sizes, it appears that the gain by optimal covariate dependent allocation probabilities increases (see Table \ref{tab:regret}). As the algorithm only gives an approximate optimal solution, we also used the approach in Section \ref{sec:lowerbound} to compute a lower bound for the limit  $n R_I$ for $n\to \infty$. The obtained lower bound of 735.1 is close to 738.2,  the corresponding regret obtained by the derived allocation function depicted in Figure \ref{fig:diets}d.

	It should be noticed that all these designs are based on the true values of the parameters and therefore require a preliminary study, or an adaptive design.

	
	
	

	\begin{table}[ht]
		\centering
		\begin{tabular}{|c|c|c|c|c|c|c|}\hline 
			& & \multicolumn{2}{c}{\% reduction}&\multicolumn{3}{|c|}{fixed allocation ($M=0$)} \\\hline 
			$n$ & $n R_I$ & $M$=4 & $M=0$&high&moderate&low\\ \hline
			$164$ &842.8&10.9 &  7.9& 0.45& 0.37& 0.18\\
			$1000$ &829.6 &13.1 &  4.4& 0.37& 0.40 & 0.23 \\
			$\infty$ & 738.2& 17.1& 3.2&0.35 &0.40& 0.25\\\hline 
		\end{tabular}	
		\caption{Diets example: $n R_I$, the percentage of reduction achieved by the optimized design with $M=0$ and $M=4$  relative to a balanced design, and the fixed allocation probabilities for $M=0$.}
		\label{tab:regret}
	\end{table}
	
	\section{Discussion\label{sec:discussion}}
	In this paper we derived optimal trial designs to determine  treatment rules that are  based on  patients' individual  covariates, assuming that the outcomes in each treatment arm  are given by a linear model. For  two treatments the optimal allocation rule depends only on the variances of the error terms in the linear models and does not depend on other model parameters, nor on the covariates. In contrast, for three or more treatments the optimal rule depends on all the linear model parameters, and the covariates.  In the latter case one faces locally optimal designs whose computation requires prior knowledge on the model parameters. If such information is not available,  extensions of the framework can be considered. For example,  sequential  designs that start with equal allocation and then compute optimal allocation rules for further stages based on interim estimates of the model parameters can be used. Another option are  Bayesian designs, where a  prior distribution for the model parameters is specified and the  regret, averaged over this prior  is minimized. Furthermore, the sequential and Bayesian approaches can be combined by computing posterior distributions of the model parameters in interim analyses and optimizing the allocation ratios based on the  regret averaged over these posterior distributions. Alternatively, a minimax approach can be used and the maximal regret across a range of plausible scenarios can be minimized. {The (approximate) optimal minimax design was computed in Section \ref{sec:approx_regret} for two treatments. The case of more than two treatment is deferred for future work.} {In the case of three treatments (or more) one can also consider the possibility
		of testing each treatment against the first, say, using the design proposed for two treatments. This would require splitting the data for each paired comparison, thus leading to loss of efficiency.}
	
	{The efficiency of a given design is measured by the regret. In Section \ref{sec:diets} we demonstrate a comparison of two designs in terms of their regrets, for which we provide exact and approximate formulas in the paper.} 
	

	A trial of the kind we study ends with recommendations on treatments, amounting to a claim of causality. Here we allow the allocation to treatments to be a function of the covariate vector $\bf X$ only, and in particular there are no other confounding variables which affects both the allocation and the response, and causality can be deduced; see, e.g., \citet{Robins} and \citet{SW}. This remains true even if the allocation functions are ``deterministic", taking only the values 0 and 1. We assume that a linear model holds true (see Section \ref{sec:compreg}) rather than taking the approach that ``all models are wrong."  This allows us to make consistent estimation of the model parameters for any allocation, barring a degenerate design which concentrates on a single $\bf X$. Assuming a model  is necessary if one wants to avoid a separate trial for each level of $\bf X$, which seems inefficient. A linear model assumption may be reasonable if the range of values of $\bf X$ is restricted, which is often very natural. As usual, linear models include a polynomial regression of some degree rather than a simple linear model (see Section \ref{sec:polynomial}), and interaction terms may be included.

	\section*{Acknowledgement}
	We thank the reviewers and the associate editor for their close reading of the manuscript and useful comments.  	 Yosef Rinott was supported in part by a grant from CIDR - Center for Interdisciplinary Data Science Research at the Hebrew University.
	Martin Posch is member of the EU Patient-Centric Clinical Trial Platform (EU-PEARL) which has received funding from the Innovative Medicines Initiative 2 Joint Undertaking, grant No 853966. This Joint Undertaking receives support from the EU Horizon 2020 Research and Innovation Program, EFPIA, Children’s Tumor Foundation, Global Alliance for TB Drug Development, and Spring Works Therapeutics. The views expressed in this publication are those of the authors. The funders and associated partners are not responsible for any use that may be made of the information contained herein.
	

		\bibliographystyle{chicago}
		\bibliography{refs}
		
		\clearpage
		
		\section{Appendix: Proofs}
		\subsubsection*{\bf Proof of Lemma \ref{lem:order}}
		We use the convexity relation $\lambda A^{-1}+(1-\lambda)B^{-1} \succeq [\lambda A+(1-\lambda)B]^{-1}$ \citep{Moore} where $A$ and $B$ are positive definite matrices. We have
		\begin{multline*}\label{eq:ineq_hetro2}
			{\bs \Sigma}_{1}+ {\bs \Sigma}_{2}={\sigma_1}
			\left(\frac{{\nu_1}{\bs Q}_{1}}{\sigma_1}\right)^{-1}+{\sigma_2}
			\left(\frac{{\nu_2}{\bs Q}_{2}}{\sigma_2}\right)^{-1} {=(\sigma_1+\sigma_2)\big[\frac{\sigma_1}{\sigma_1+\sigma_2}
				\left(\frac{{\nu_1}{\bs Q}_{1}}{\sigma_1}\right)^{-1}+\frac{\sigma_2}{\sigma_1+\sigma_2}
				\left(\frac{{\nu_2}{\bs Q}_{2}}{\sigma_2}\right)^{-1}\big]}\\
			\succeq
			(\sigma_1+\sigma_2)\Big[\frac{\sigma_1}{\sigma_1+\sigma_2}
			\left(\frac{{\nu_1}{\bs Q}_{1}}{\sigma_1}\right)+\frac{\sigma_2}{\sigma_1+\sigma_2}\left(\frac{{\nu_2}{\bs Q}_{2}}{\sigma_2}\right)\Big]^{-1}
			=
			(\sigma_1+\sigma_2)^2 \mathbb{\bs Q}^{-1},
		\end{multline*}
		where in the last equality we used the fact that
		$\nu_1{\bs Q}_{1}+{\nu_2}{\bs Q}_{2}=
		\mathbb{\bs Q}$.
		Equality holds when the two matrices to the {right} of the $\succeq$ sign are  equal, which happens when
		$\pi_1({\bf x})=\frac{\sigma_1}{\sigma_1+\sigma_2}$ and $\pi_2({\bf x})=\frac{\sigma_2}{\sigma_1+\sigma_2}$. \qed

		\subsubsection*{Proof of Theorem \ref{thm:mult_K_2}} 
		We start with the case of $p=1$.     
		It is enough to show (as the other part of the integral, from $\theta$ to 1, is symmetric) that
		\begin{equation}\label{eq:part_int}
			n^{3/2-\ve}  \int_0^\theta \left| P(\widehat{g}_2(x) > \widehat{g}_1(x))-\Phi\left(\frac{-\sqrt{n}[g_1(x) - g_2(x)]}{\sqrt{V(x)}}\right)\right|[g_1(x) - g_2(x)]   f(x)dx \to 0.
		\end{equation}
		The idea of the proof is quite simple. For $x \in \theta \pm 1/\sqrt{n}$ the quantity in absolute value converges to zero {at rate $1/\sqrt{n}$}, and both the range of the integral and $g_1(x) - g_2(x)$ are of order $1/\sqrt{n}$ each.
		For other $x$'s both the probability and $\Phi$ in the integrand are exponentially small in $n$.

		We will prove \eqref{eq:part_int} by showing that 
		\begin{equation}\label{eq:part_int_1}
			n^{3/2-\ve}  \int_0^\theta \left| P(\widehat{g}_2(x) > \widehat{g}_1(x))-E\Phi\left(\frac{-\sqrt{n}(\beta_2-\beta_1)(\theta-x)}{S(x)}\right)\right|[g_1(x) - g_2(x)]   f(x)dx \to 0\text{ and}
		\end{equation}
		\begin{equation}\label{eq:part_int_2}
			n^{3/2-\ve}  \int_0^\theta \left|E\Phi\left(\frac{-\sqrt{n}(\beta_2-\beta_1)(\theta-x)}{S(x)}\right)-\Phi\left(\frac{-\sqrt{n}[g_1(x) - g_2(x)]}{\sqrt{V(x)}}\right)\right|[g_1(x) - g_2(x)]   f(x)dx \to 0,
		\end{equation}
		where
		\begin{equation}\label{eq:S^2}
			S^2(x)=n(1,x)\left\{ \sigma_1^2 \left( \sum_{i \in T_1} \binom{1}{X_i} (1,X_i)\right)^{-1}    + \sigma_2^2 \left( \sum_{i \in T_2} \binom{1}{X_i} (1,X_i)\right)^{-1} \right\}\binom{1}{x},
		\end{equation}
		and the expectation in \eqref{eq:part_int_1} and \eqref{eq:part_int_2} is with respect to $X_1,\ldots,X_n$ which appear in $S^2(x)$ and the random sets $T_k$.
		We start with \eqref{eq:part_int_2}. 
		Since $\theta$ is the intersection point, $\alpha_1+\beta_1 \theta = \alpha_2+\beta_2 \theta$; therefore,
			$g_1(x)-g_2(x)=\alpha_1+\beta_1 x- (\alpha_2 + \beta_2 x)=(\beta_2 - \beta_1) (\theta-x)$.
			Hence, \eqref{eq:part_int_2} can be written as
			\begin{equation*}
				n^{3/2-\ve}  \int_0^\theta \left| E\Phi\left(\frac{-\sqrt{n}(\beta_2-\beta_1)(\theta-x)}{ S(x)}\right)-\Phi\left(\frac{-\sqrt{n}(\beta_2-\beta_1)(\theta-x)}{\sqrt{ V(x)}}\right)\right|(\beta_2-\beta_1)(\theta-x)   f(x)dx \to 0.
			\end{equation*}
			Fix $c_n=n^{\ve/5}$. We will show that
			\begin{equation}\label{eq:part_int1}
				n^{3/2-\ve}  \int_0^{\theta-c_n/\sqrt{n}} \left| E\Phi\left(\frac{-\sqrt{n}(\beta_2-\beta_1)(\theta-x)}{S(x)}\right)-\Phi\left(\frac{-\sqrt{n}(\beta_2-\beta_1)(\theta-x)}{\sqrt{ V(x)}}\right)\right|(\beta_2-\beta_1)(\theta-x)   f(x)dx \to 0,
			\end{equation}
			and that
			\begin{equation}\label{eq:part_int2}
				n^{3/2-\ve}  \int_{\theta-c_n/\sqrt{n}}^\theta \left| E\Phi\left(\frac{-\sqrt{n}(\beta_2-\beta_1)(\theta-x)}{S(x)}\right)-\Phi\left(\frac{-\sqrt{n}(\beta_2-\beta_1)(\theta-x)}{\sqrt{ V(x)}}\right)\right|(\beta_2-\beta_1)(\theta-x)   f(x)dx \to 0.
			\end{equation}
			We show \eqref{eq:part_int1} by arguing that the probabilities in \eqref{eq:part_int1} are exponentially small. 
			With notation defined in Section \ref{sec:compreg} we have that
			\[
			V(x) \le (1+x^2) \left(\frac{\sigma_1^2}{\nu_1 \lambda_{min}({\bs Q_1}) }+
			\frac{\sigma_2^2}{\nu_2 \lambda_{min}({\bs Q_2}) }\right) \le C,
			\]
			for a constant $C>0$ (since $x$ is bounded).
			With the bound $\Phi(-t) \le \exp(-t^2/2)$ for  $t>1$  and
			$(\theta-x)\sqrt{n} \le c_n$ we have
			\[
			\Phi\left(\frac{-(\beta_2-\beta_1)(\theta-x)\sqrt{n}}{\sqrt{ V(x)}}\right) \le \exp\left(-\frac{(\beta_2-\beta_1)^2c_n^2}{2C}\right),
			\]
			which shows that the second normal probability $\Phi(\cdot)$ in \eqref{eq:part_int1} is exponentially small. We now argue that the expectation in \eqref{eq:part_int1} is also exponentially small. Consider the event 
			\begin{multline} \label{eq:A_n}
				A_n:=\left\{ \lambda_{min}\left(\frac{1}{n} \sum_{i \in T_1} \binom{1}{X_i} (1,X_i)\right)  \le \nu_1 \frac{\lambda_{min}({\bs Q_1})}{2}\right\}\\
				\bigcup \left\{ \lambda_{min}\left(\frac{1}{n} \sum_{i \in T_2} \binom{1}{X_i} (1,X_i)\right)  \le \nu_2 \frac{\lambda_{min}({\bs Q_2})}{2} \right\}.    
			\end{multline}
			Writing $\frac{1}{n} \sum_{i \in T_1} \binom{1}{X_i} (1,X_i)=\frac{1}{n} \sum_{i=1}^n \binom{1}{X_i} (1,X_i)I(i \in T_1)$,   a sum of iid terms, we can apply Theorem 5.1 of \citet{Tropp} and the union bound to conclude that $P(A_n) \le C\exp(-Cn)$ for some constant $C>0$. Going back to \eqref{eq:part_int1}, write
			\begin{multline*}
				E\Phi\left(\frac{-\sqrt{n}(\beta_2-\beta_1)(\theta-x)}{S(x)}\right) \\
				=
				E\Phi\left(\frac{-\sqrt{n}(\beta_2-\beta_1)(\theta-x)}{S(x)}\right)I(A_n) + E\Phi\left(\frac{-\sqrt{n}(\beta_2-\beta_1)(\theta-x)}{S(x)}\right)I(A_n^c)\\
				\le P(A_n) + E\Phi\left(\frac{-\sqrt{n}(\beta_2-\beta_1)(\theta-x)}{S(x)}\right)I(A_n^c).    
			\end{multline*}
			When $A_n^c$ occurs, the denominator $S^2(x)$ is bounded and the same argument as above shows that this probability is exponentially small. We conclude that both probabilities in \eqref{eq:part_int1} are exponentially small, which implies that \eqref{eq:part_int1} is correct.
			
			We now prove \eqref{eq:part_int2}. By the mean value theorem, for $x \in (\theta-c_n/\sqrt{n},\,\theta)$ and for some $V^*(x)$ between $S^2(x)$ and $V(x)$, we have for some $C>0$
			\begin{multline*}
				\left|\Phi\left(\frac{-\sqrt{n}(\beta_2-\beta_1)(\theta-x)}{\sqrt{S^2(x)}}\right)-\Phi\left(\frac{-\sqrt{n}(\beta_2-\beta_1)(\theta-x)}{\sqrt{ V(x)}}\right)\right|\\
				=|S^2(x)-V(x)|\varphi\left(\frac{-\sqrt{n}(\beta_2-\beta_1)(\theta-x)}{\sqrt{ V^*(x)}}\right)\frac{\sqrt{n}(\beta_2-\beta_1)(\theta-x)}{2} (V^*(x))^{-3/2} \\
				\le {C}|S^2(x)-V(x)|c_n (\min\{S^2(x),V(x)\})^{-3/2},\text{ and} 
			\end{multline*}
			\begin{equation}\label{eq:M_1_M_2}
				\{S^2(x)-V(x)\}= (1,x)\left\{ \frac{\sigma_1^2}{\nu_1}\left( {\bf M}_1^{-1}  - {\bs Q}_{1}^{-1}\right)  + \frac{\sigma_2^2}{\nu_2} \left({\bf M}_2^{-1}-{\bs Q}_{2}^{-1}\right) \right\}\binom{1}{x},
			\end{equation}
			where 
			${\bf M}_k=\frac{1}{n\nu_k} \sum_{i \in T_k} \binom{1}{X_i} (1,X_i)$, $k=1,2$.
			We have that 
			${\bf M}_1^{-1}  - {\bs Q}_{1}^{-1}={\bf M}_1^{-1}\left({\bs Q}_{1}-{\bf M}_1\right){\bs Q}_{1}^{-1}$,
			implying that the difference of the inverses is small when the difference of the matrices is small and the inverse is bounded.  
			Define the event
			\[
			\widetilde{A}_n := \{ \max\{ \| {\bs Q}_{1}-{\bf M}_1 \|_\infty, \| {\bs Q}_{2}-{\bf M}_2 \|_\infty\} \ge c_n/\sqrt{n} \},
			\]
			where $\| \cdot \|_\infty$ denotes the maximum (entry-wise) norm. By Hoeffding's inequality (since the random variables are bounded), $P(\widetilde{A}_n) \le C\exp(-C c_n^2)$ (which is exponentially small). Therefore, as before, this event can be ignored. On the complement of $\widetilde{A}_n$, ${\bf M}_1^{-1}$, ${\bf M}_2^{-1}$, $S^2(x)$ are all bounded, and by similar arguments $S^2(x)$ and $V(x)$ are bounded below. Hence, when $\widetilde{A}_n^c$ occurs,
			\[
			\left|\Phi\left(\frac{-\sqrt{n}(\beta_2-\beta_1)(\theta-x)}{S(x)}\right)-\Phi\left(\frac{-\sqrt{n}(\beta_2-\beta_1)(\theta-x)}{\sqrt{ V(x)}}\right)\right| \le C \frac{c_n^2}{\sqrt{n}}. 
			\]
			Therefore the integral in \eqref{eq:part_int2} is bounded by
			\[
			C {n}^{1-\ve} c_n^2 \int_{\theta-c_n/\sqrt{n}}^\theta (\theta-x)   f(x)dx \le C n^{1/2-\ve} c_n^3 \int_{\theta-c_n/\sqrt{n}}^\theta   f(x)dx  \le C c_n^4/{n^{\ve}} =n^{-\ve/5} \to 0,
			\]
			which completes the proof of \eqref{eq:part_int2} and hence of \eqref{eq:part_int_2}.
			
			We now show \eqref{eq:part_int_1} by 
			dividing the integral into two ranges: \textbf{(a)}  $[\theta-c_n/\sqrt{n}, \theta]$, and \textbf{(b)}  $[0, \theta-c_n/\sqrt{n}]$.

			
			We have that 
			\begin{multline*}
				P\left(\widehat{g}_2(x)>\widehat{g}_1(x)\right)=P\left( \left[ \binom{\widehat{\alpha}_1 -\alpha_1}{\widehat{\beta}_1-\beta_1}-\binom{\widehat{\alpha}_2 -\alpha_2}{\widehat{\beta}_2-\beta_2} \right]^t \binom{1}{x} > (\beta_2-\beta_1)(x-\theta)\right)\\
				=P\left( \left[ \frac{{\bf M}_1^{-1}}{n\nu_1} \sum_{i \in T_1} \binom{\ve_i}{X_i \ve_i} - \frac{{\bf M}_2^{-1}}{n\nu_2} \sum_{i \in T_2} \binom{\ve_i}{X_i \ve_i} \right]^t \binom{1}{x} > (\beta_2 - \beta_1)(x-\theta) \right)\\
				=P\left(\sum_{i=1}^n a_i \ve_i > (\beta_2-\beta_1)(x-\theta)\right),
			\end{multline*}
			where $a_i = \left\{ \begin{array}{cc} (1,X_i)\frac{{\bf M}_1^{-1}}{n\nu_1} \binom{1}{x} & i \in T_1\\
				-(1,X_i)\frac{{\bf M}_2^{-1}}{n\nu_2} \binom{1}{x} & i \in T_2 \end{array} \right.$,\,
			${\bf M}_1,{\bf M}_2$ are defined in \eqref{eq:M_1_M_2} 	and $\ve_i$ are the regression error terms.  Notice that $Var(\ve_i)$ is equal to $\sigma_1^2$ if $i \in T_1$ and to $\sigma_2^2$ if $i \in T_2$. We are going to condition on $X_1,\ldots, X_n$ obtained in the experiment and their allocation, denoted together by $D$. Then
			$Var(\sum_{i=1}^n a_i \ve_i | D) = S^2(x)/n$, where $S^2(x)$ is defined in \eqref{eq:S^2}. 
			
			We aim to apply the Berry Esseen Theorem, for which we need some bounds. The values $na_i$ are bounded below since the matrices ${\bf M}_k$ are bounded above.
			In order to obtain upper bounds, recall the event $A_n$  defined in \eqref{eq:A_n}. On the set $A_n^c$ we have that $na_i$ are bounded above. Since $P(A_n)$ is exponentially small, we can ignore it, and assume that $na_i$ are bounded with probability one. We can apply 
			the Berry-Esseen Theorem for the non-identically distributed  case (for a convenient reference see \citet{CGS} (3.27)) to  $\sum_{i=1}^n \frac{\sqrt{n}a_i}{S(x)} \varepsilon_i \equiv \sum_{i=1}^n V_i$ where $\sum_{i=1}^n Var(V_i | D)=1$ and $E(|V_i|^3 \mid D) \le max\{na_i/S(x) E(|\varepsilon_i|^3 | D): 1 \le i \le n\}n^{-3/2}$ to obtain by unconditioning on $D$
			\[
			\left|P\left(  \sum_{i=1}^n a_i \ve_i > (\beta_2-\beta_1)(x-\theta) \right) -E \Phi\left(\frac{-\sqrt{n}(\beta_2-\beta_1)(\theta-x)}{S(x)}\right)\right| \le C/\sqrt{n}.
			\]
			Therefore,
			\begin{multline}\label{eq: bezalel}
				n^{3/2-\ve}  \int_{\theta-c_n/\sqrt{n}}^\theta \left| P(\widehat{g}_2(x)>\widehat{g}_1(x))-E\Phi\left(\frac{-\sqrt{n}(\beta_2-\beta_1)(\theta-x)}{S(x)}\right)\right|(\beta_2-\beta_1)(\theta-x)f(x)dx
				\\
				\le n^{3/2-\ve} \frac{C}{\sqrt{n}} \int_{\theta-c_n/\sqrt{n}}^\theta (\beta_2-\beta_1)(\theta-x)f(x)dx\le  n^{3/2-\ve} \frac{C}{\sqrt{n}}  \cdot \frac{c_n}{\sqrt{n}} \int_{\theta-c_n/\sqrt{n}}^\theta f(x)dx\\
				\le C \frac{c_n^2}{{n}^\ve}= C \frac{n^{2\ve/5}}{{n}^\ve} \to 0.
			\end{multline}
			By a standard large deviation bound, the  integral of \eqref{eq:part_int_1}  from 0 to $\theta-c_n/\sqrt{n}$ multiplied by $n^{3/2-\ve}$ is easily shown to be of order $n^{3/2-\ve} O(exp(-c_n \sqrt{n})) \to 0$. This proves the theorem for covariates in $\mathbb{R}$ ($p$=1).
			
			We now consider the case of covariates in $\mathbb{R}^p$ for $p\ge 2$.
			Suppose, without loss of generality, that ${\beta}_{2,1}>{\beta}_{1,1}$. For given ${\bf x}_{-1}=(x_2,\ldots,x_p)$ define $\theta_1({\bf x}_{-1}):=\frac{\alpha_1-\alpha_2+({\bs \beta}_{1,-1}-{\bs \beta}_{2,-1})^t {\bf x}_{-1} }{\beta_{2,1}-\beta_{1,1}}$. For ${\bf x}$ such that $\theta_1({\bf x}_{-1}) \notin [0,1]$ the probability of making a mistake (i.e., $\widehat{g}_2({\bf x}) - \widehat{g}_1({\bf x})$ has the wrong sign) is exponentially small. Therefore,
			\begin{multline*}
				R(\pi_1,\pi_2)= \int_{[0,1]^{p-1}}\int_0^{\theta_1} P(\widehat{g}_2(\x) > \widehat{g}_1(\x))(g_1(\x) - g_2(\x))f(\x)dx_1 I(\theta_1 \in [0,1])d{\x}_{-1}\\+
				\int_{[0,1]^{p-1}}\int_{\theta_1}^1 P(\widehat{g}_2(\x) > \widehat{g}_1(\x))(g_1(\x) - g_2(\x))f(\x)dx_1 I(\theta_1 \in [0,1])d{\x}_{-1}+a_n,  
			\end{multline*}
			where $a_n$ is exponentially small. By a slight variation on the above one-dimensional case, the inner integral equals 
			\begin{equation}\label{eq:inner_int}
				\int_0^1 \Phi\left(\frac{-\sqrt{n}|g_2({\bf x}) - g_1({\bf x})|}{\sqrt{V({\bf x})}}\right)|g_2({\bf x}) - g_1({\bf x})| f({\bf x})dx_1I(\theta_1 \in [0,1])+o(1/n^{3/2-\ve}),
			\end{equation}
			where the error term $o$ is uniform in  ${\bf x}_{-1}$.
			Taking the outer integral, the theorem follows. \qed

			\subsubsection*{Proof of Theorem \ref{thm:one_dim_lim}}
			As computed two lines below  \eqref{eq:S^2}, $g_1(x)-g_2(x)=(\beta_2-\beta_1)(\theta-x)$. We will show that for continuous $t$
			\begin{equation}\label{eq:lim_t_infty}
				\lim_{t \to \infty} t \int_0^\theta \Phi\left(\frac{-\sqrt{t}(\beta_2-\beta_1)(\theta-x)}{\sqrt{ V(x)}}\right)(\beta_2-\beta_1)(\theta-x)   f(x)dx = \frac{ V(\theta)f(\theta)}{4 (\beta_2-\beta_1)}. 
			\end{equation}
			The integral from $\theta$ to 1 is similar, yielding the desired limit.
			By L'H\^opital's rule the limit in \eqref{eq:lim_t_infty} equals
			\[
			\lim_{t \to \infty} \frac{t^2}{2} \int_0^\theta\varphi\left(\frac{\sqrt{t}(\beta_2-\beta_1)(\theta-x)}{\sqrt{ V(x)}}\right) \frac{(\beta_2-\beta_1)(\theta-x)}{\sqrt{ V(x)}}\frac{1}{\sqrt{t}}(\beta_2-\beta_1)(\theta-x) f(x) dx.
			\]
			Substitution of $x$ by $y=(\theta-x)\sqrt{t}$ in the integral yields 
			\[
			\lim_{t \to \infty} \frac{1}{2} \int_0^{\theta\sqrt{t}}\varphi\left(\frac{(\beta_2-\beta_1)y}{\sqrt{ V(\theta-y/\sqrt{t})}}\right) \frac{(\beta_2-\beta_1)^2}{\sqrt{ V(\theta-y/\sqrt{t})}}y^2 f(\theta-y/\sqrt{t}) dy.
			\]
			The limit is
			\[
			\frac{1}{2} \int_0^{\infty}\varphi\left(\frac{(\beta_2-\beta_1)y}{\sqrt{ V(\theta)}}\right) \frac{(\beta_2-\beta_1)^2}{\sqrt{V(\theta)}}y^2 f(\theta) dy.
			\]
			Substitution of $y$ with $z=\frac{(\beta_2-\beta_1)y}{\sqrt{ V(\theta)}}$ in the integral yields the claimed limit
			$\frac{f(\theta) V(\theta)}{2(\beta_2-\beta_1)} \int_0^{\infty}\varphi\left(z\right) z^2  dz=\frac{f(\theta)  V(\theta)}{4(\beta_2-\beta_1)}.$ \qed

			\subsubsection*{Proof of Theorem \ref{thm:2}}
			Consider the inner integral in \eqref{eq:inner_int}. By Theorem \ref{thm:one_dim_lim}, the limit of the inner integral times $n$ is
			\[
			\frac{V(\theta_1,{\bf x}_{-1}) f(\theta_1,{\bf x}_{-1})}{2(\beta_{2,1}-\beta_{1,1})}. 
			\]
			By taking the outer integral the result follows. \qed
			
			\subsection*{The limit of the regret for polynomial regression}
			Suppose that $g_1(x)=\alpha_1 + (x,\ldots,x^J)^t{\bs \beta}_1$ and $g_2(x)=\alpha_2+(x,\ldots,x^J)^t{\bs \beta}_2$. Let $\theta$ be an intersection point, i.e., $g_1(\theta)=g_2(\theta)$. We have 
			\[
			g_2(x)-g_1(x)=g_2(x)-g_2(\theta)-[g_1(x)-g_1(\theta)]={\bs \eta}^t({\bs \beta_2}-{\bs \beta_1}),
			\]
			where ${\bs \eta}:=(x-\theta,x^2-\theta^2,\ldots,x^J-\theta^J)^t$. The following proposition is parallel to Theorem \ref{thm:one_dim_lim} (in the one-dimensional case) and provides the limit of the regret; see \eqref{eq:polya}. The proof is similar, although some modifications are required.
			\begin{proposition}\label{prop:poly}
				For a continuous argument $s$
				\begin{equation}\label{eq:lim_t_infty_poly}
					\lim_{s \to \infty} s \int_0^\theta \Phi\left(\frac{-\sqrt{t}({\bs \beta_2}-{\bs \beta_1})^t{\bs \eta}}{\sqrt{ V(x)}}\right)({\bs \beta_2}-{\bs \beta_1})^t{\bs \eta} f(x)dx = \frac{f(\theta)V(\theta)}{4|({\bs \beta_2}-{\bs \beta_1})^t{\bs \zeta}|}, 
				\end{equation}
				where ${\bs \zeta}:=(1,2\theta,...,J\theta^{J-1})^t$.
			\end{proposition}
			{\bf Proof.} By L'H\^opital's rule the limit in \eqref{eq:lim_t_infty_poly} equals
			\[
			\lim_{s \to \infty} \frac{s^2}{2} \int_0^\theta\varphi\left(\frac{\sqrt{t}({\bs \beta_2}-{\bs \beta_1})^t{\bs \eta}}{\sqrt{ V(x)}}\right) \frac{({\bs \beta_2}-{\bs \beta_1})^t{\bs \eta}}{\sqrt{ V(x)}}\frac{1}{\sqrt{t}}({\bs \beta_2}-{\bs \beta_1})^t{\bs \eta} f(x) dx.
			\]
			Substituting $y=(\theta-x)\sqrt{s}$ we get 
			\[
			\sqrt{s}{\bs \eta}=\sqrt{s}\begin{pmatrix} x-\theta \\ x^2-\theta^2 \\ \vdots \\ x^J-\theta^J \end{pmatrix}=\sqrt{s}\begin{pmatrix} x-\theta \\ (x-\theta)(x+\theta)\\ \vdots \\ (x-\theta)(x^{J-1}+x^{J-2}\theta+\ldots+\theta^{J-1}) \end{pmatrix}=y\widetilde{\bs \eta},
			\]
			where 
			\begin{multline*} 
				\widetilde{\bs \eta}:=(1,x+\theta,\ldots,x^{J-1}+x^{J-2}\theta+\ldots+\theta^{J-1})^t\\
				=(1,2\theta-y/\sqrt{s},\ldots,(\theta-y/\sqrt{s})^{J-1}+(\theta-y/\sqrt{s})^{J-2}\theta+\ldots+\theta^{J-1})^t.
			\end{multline*}
			The integral reads
			\[
			\lim_{s \to \infty} \frac{1}{2} \int_0^{\theta\sqrt{s}}\varphi\left(y\frac{({\bs \beta_2}-{\bs \beta_1})^t\widetilde{\bs \eta}}{\sqrt{ V(\theta-y/\sqrt{s})}}\right) \frac{y^2[({\bs \beta_2}-{\bs \beta_1})^t\widetilde{\bs \eta}]^2}{\sqrt{ V(\theta-y/\sqrt{s})}} f(\theta-y/\sqrt{t}) dy.
			\]
			The limit is
			\[
			\frac{1}{2} \int_0^{\infty}\varphi\left(\frac{y({\bs \beta_2}-{\bs \beta_1})^t{\bs \zeta}}{\sqrt{ V(\theta)}}\right) \frac{[y({\bs \beta_2}-{\bs \beta_1})^t{\bs \zeta}]^2}{\sqrt{ V(\theta)}} f(\theta) dy.
			\]
			Substitution of $y$ with $z=\frac{y({\bs \beta_2}-{\bs \beta_1})^t{\bs \zeta}}{\sqrt{ V(\theta)}}$ yields
			\[
			\frac{f(\theta) V(\theta)}{2|({\bs \beta_2}-{\bs \beta_1})^t{\bs \zeta}|} \int_0^{\infty}\varphi\left(z\right) z^2  dz=\frac{f(\theta)V(\theta)}{4|({\bs \beta_2}-{\bs \beta_1})^t{\bs \zeta}|},
			\]
			which completes the proof. \qed

			For the proofs of Theorems \ref{thm:K_one_dim1} and \ref{thm:K_one_dim_lim} we need the lemma below and some notation. 
			For every $m$ and $k$, let 
			$I_{m,k}:= \int_{\theta_{m-1}}^{\theta_{m}}P\left(\widehat{g}_k(x) > \max_{\ell \ne k} \widehat{g}_\ell(x) \right)\left[ g_{m}({x}) - g_k({x})\right]f({x})d{x}$. We then have $R(\pi_1,\ldots,\pi_K)=\sum_{k=1}^K\sum_{m=1}^K I_{m,k}$.

			\begin{lemma}
				\label{thm:K_one_dim}
				Under the assumptions in the beginning of Section \ref{sec:Knorm}, the integral $I_{m,k}$
				is exponentially small for $k \ne m-1,m+1$, and for every $\ve>0$
				\begin{eqnarray}
					&\lim_{n \to \infty} n^{3/2-\ve} \left| I_{m,m-1}-\int_{\theta_{m-1}}^{\theta_{m-1}+\frac{c_n}{\sqrt{n}}}\Phi\left(\frac{-\sqrt{n}\{g_m(x) - g_{m-1}(x)\}}{\sqrt{V_{m-1}(x)}}\right)\left[ g_{m}({x}) - g_{m-1}({x})\right]f({x})d{x}\right|=0, \label{eq:m,m-1} \\\nonumber
					&\lim_{n \to \infty} n^{3/2-\ve} \left| I_{m,m+1}-\int_{\theta_{m}-\frac{c_n}{\sqrt{n}}}^{\theta_{m}}\Phi\left(\frac{-\sqrt{n}\{g_m(x) - g_{m+1}(x)\}}{\sqrt{V_{m}(x)}}\right)\left[ g_{m}({x}) - g_{m+1}({x})\right]f({x})d{x}\right|=0,
				\end{eqnarray}
				where $c_n = n^{\ve/5}$ and $V_m(x)=(1,x)\left( {\bs \Sigma}_{m}+ {\bs \Sigma}_{m+1} \right)\binom{1}{x}$. Moreover, for every $m$ and $k$
				\[
				n^{3/2-\ve} \left| I_{m,k}
				-\int_{\theta_{m-1}}^{\theta_{m}}\left(\int \prod_{\ell=1,\ldots,K,l\not=k} \Phi\left(\frac{z\xi_k(x)+\sqrt{n}[g_k({x})-g_\ell({x})]}{\xi_\ell(x)}\right)\varphi(z)dz\right)\left[ g_{m}({x}) - g_k({x})\right]f({x})d{x}\right| \to 0,
				\]
				as $n \to \infty$.
			\end{lemma}

			\subsubsection*{\bf Proof of Lemma \ref{thm:K_one_dim}}
			
			
			The proof of the first claim in Lemma \ref{thm:K_one_dim} is similar to that of Theorem \ref{thm:mult_K_2}. Here is a sketch. For $x \in I_{m,k}$ and bounded away from $\theta_{m-1}$ and $\theta_m$, and $k\ne m$
			both the quantities
			\[
			P\left(\widehat{g}_k(x) > \max_{\ell \ne k} \widehat{g}_\ell(x) \right)\text{ and }
			\int \prod_{\ell=1,\ldots,K,l\not=k} \Phi\left(\frac{z\xi_k(x)+\sqrt{n}[g_k({x})-g_\ell({x})]}{\xi_\ell(x)}\right)\varphi(z)dz
			\]
			are exponentially small (for $k=m$ the regret is zero). The same holds true for every $x$ and when $k\ne m-1,m+1$. 
			
			We now prove \eqref{eq:m,m-1} concerning $I_{m,m-1}$ (where $m>1$). The other relation for 
			$I_{m,m+1}$ is similar.
			Fix $c_n=n^{\ve/5}$ and consider $x$ such that $x\in (\theta_{m-1},\theta_{m-1}+c_n/\sqrt{n})$ . For such $x$ we have by a standard large deviations argument that
			\begin{multline}\label{eq:b_n}
				P\left(\widehat{g}_{m-1}(x) > \max_{\ell \ne {m-1}} \widehat{g}_\ell(x) \right) =  P(\widehat{g}_{m-1}(x)> \widehat{g}_{m}(x)) + a_n, \text{ and }\\
				\int \prod_{\ell=1,\ldots,K,l\not=m-1} \Phi\left(\frac{z\xi_{m-1}(x)+\sqrt{n}[g_{m-1}({x})-g_\ell({x})]}{\xi_\ell(x)}\right)\varphi(z)dz  =\Phi\left(\frac{-\sqrt{n}[g_m(x) - g_{m-1}(x)]}{\sqrt{V_{m-1}(x)}}\right)+b_n,
			\end{multline}
			where $a_n$ and $b_n$ are exponentially small (uniformly in $x$). By a Berry-Esseen type bound applied to the first part of the integrand, and the smallness of the other part for $x$ is the given range of the integral, and the smallness of the range itself , as in \eqref{eq: bezalel}, we obtain
			\[
			\int_{\theta_{m-1}}^{\theta_{m-1}+\frac{c_n}{\sqrt{n}}}\left[P(\widehat{g}_{m-1}(x)>\widehat{g}_{m}(x))- \Phi\left(\frac{-\sqrt{n}[g_m(x) - g_{m-1}(x)]}{\sqrt{V_m(x)}}\right)\right]\left[ g_{m}({x}) - g_{m-1}({x})\right]f({x})d{x}=o(1/n^{3/2-\ve}),
			\]
			and \eqref{eq:m,m-1} follows. This, together with \eqref{eq:b_n}, implies the last part of the theorem  when $k=m-1$.  \qed
			
			\subsubsection*{\bf Proof of Theorem \ref{thm:K_one_dim1}} This follows directly from the last part of Lemma \ref{thm:K_one_dim} and  the fact that $R(\pi_1,\ldots,\pi_K)=\sum_{k,m =1}^K I_{m,k}$.  Note that we can replace that latter sum by $\sum_{m =1}^K I_{m,m-1}$
			since $I_{m,k}$
			is exponentially small for $k \ne m-1,m+1$.

			\subsection*{Proof of Theorem \ref{thm:K_one_dim_lim}}
			The first part of Lemma \ref{thm:K_one_dim} implies that
			\[
			\sum_{k=1}^K\sum_{m=1}^K I_{m,k}=\sum_{m=1}^{K-1}\int_{\theta_{m}-\frac{c_n}{\sqrt{n}}}^{\theta_{m}+\frac{c_n}{\sqrt{n}}}\Phi\left(\frac{-\sqrt{n}|g_m(x) - g_{m-1}(x)|}{\sqrt{V_{m-1}(x)}}\right)\left| g_{m}({x}) - g_{m-1}({x})\right|f({x})d{x}+o(1/n).
			\]  
			The limit of the latter integrals can be calculated using Theorem \ref{thm:one_dim_lim} implying the result. \qed


	\end{document}